\documentclass{amsart}

\usepackage{graphics,verbatim,color,amsthm}
\usepackage[all]{xy}
\usepackage{MnSymbol}
\usepackage{epsfig}% Added by Andreas
\usepackage{hyperref}% Added by Andreas
\hypersetup{colorlinks,pdfauthor={JKM}}% Added by Andreas

\theoremstyle{plain}
\newtheorem{theorem}{Theorem}
\newtheorem{backtheorem}{Background Theorem}

\newtheorem{proposition}{Proposition}

\newtheorem{lemma}{Lemma}  
\newtheorem{cor}{Corollary} 
\newtheorem{definition}{Definition}

\newtheorem{acknowledgement}{Acknowledgement}

\theoremstyle{remark}
\newtheorem*{remark}{Remark}

\theoremstyle{definition}
\newtheorem{example}{Example}

\def\R{\mathbb{R}}
\def\E{\mathbb{E}}

\def\C{\mathbb{C}}

\def\bP{\mathbb{P}}

\def\mL{{\mathcal L}}
\def\mB{{\mathcal B}}
\def\mR{{\mathcal R}}

 \newcommand{\dd}[2]
{
{{\partial #1}   \over {\partial #2}}
}
\makeatletter
\g@addto@macro{\endabstract}{\@setabstract}
\newcommand{\authorfootnotes}{\renewcommand\thefootnote{\@fnsymbol\c@footnote}}%
\makeatother

\definecolor{gray}{rgb}{.8,.8,.8}
\definecolor{comment}{rgb}{0,0.5,.8}

\newcommand{\LINES}{{\rm\it LINES}}

\newcommand {\eh}{{\textstyle \frac{1}{2}}}

\begin{document}

\title{ Lagrangian Relations and Linear  billiards} 

\maketitle

\begin{center}
  \large
  \authorfootnotes
 Jacques F{\'e}joz\textsuperscript{1}, 
  Andreas Knauf\textsuperscript{2},
  Richard Montgomery\textsuperscript{3}\par \bigskip
 \normalsize
 \textsuperscript{1}Universit{\'e} Paris-Dauphine \& Observatoire de Paris,
 \texttt{jacques.fejoz@dauphine.fr} \par
 \textsuperscript{2}Department of Mathematics,
Friedrich-Alexander-University Erlangen-N{\"u}rnberg,
Cauerstr.\ 11, D-91058 Erlangen, 
Germany, \texttt{knauf@math.fau.de}
 \par
\textsuperscript{3}Mathematics Department,
UC Santa Cruz, 
4111 McHenry, 
Santa Cruz, CA 95064, USA, \texttt{rmont@ucsc.edu}  \bigskip

  \today
\end{center}

%CONNOR corrected some typos...

 \begin{abstract} 

  Motivated by the high-energy limit of
the $N$-body problem we   construct  non-deterministic billiard process.  The billiard  table is the complement of a   finite collection 
 of linear subspaces  within a  Euclidean vector space.   A trajectory  is a    constant speed polygonal curve with vertices
 on the subspaces and  change  of   direction  upon hitting   a subspace   governed by   ``conservation of momentum''
(mirror reflection).  The  itinerary of a trajectory  is the list of subspaces it hits, in   order.   
%Two basic questions are: 
 (A)  Are itineraries  finite?
 (B)  What is the structure  of the space of all   trajectories having   a   fixed itinerary?
 In a beautiful series of papers Burago-Ferleger-Kononenko [BFK] answered (A) affirmatively by using non-smooth metric geometry ideas and the notion of a Hadamard space. 
 We answer (B) by proving that this  space of   trajectories  is diffeomorphic to   a Lagrangian relation on the space of lines in the Euclidean  space.
Our methods combine those of BFK with the notion of a   generating family for a   Lagrangian relation. 
 \end{abstract} 

\section{ Introduction.} \label{sect:1}

. 

\subsection{Euclidean Data. Point Billiards.  Motivating Example.}

%{\coma{[hi]}}

Consider a Euclidean vector space  $E$  endowed 
 with a   finite collection ${\mathcal L}$
 of linear 
 subspaces which we call ``collision subspaces''.  
 %Assume the subspaces all have the same  {\coma why?} codimension $d$, $d \ge 1$.  
 Write 
   \begin{equation}
 \label{collision}
 C =   \bigcup_{ L \in {\mathcal L}}   L   \qquad \text{ (collision locus)}
 \end{equation}
for the collision locus and
\begin{equation}
 \label{table}
 E^0 = E \setminus C  \qquad \text{ (our billiard  table)} 
 \end{equation}
 for its complement.  Play   billiards on  $E^0$ !  
 
A  ``billiard trajectory''  will be  a  certain type of polygonal curve  $q: \R \to \E$ all of
whose vertices are collisions, i.e.   lie  on $C$. 
%We parameterize $q$ at constant speed.
 When $q$ hits a subspace $L \in {\mathcal L}$ it  switches directions  by  bouncing off  of  $L$ according
 to  the laws of  reflection  (see equations  (\ref{eq:energy}), (\ref{eq:momentum}) below). 
 %Thus the $L$'s are codimension $d$   mirrors.
Imagine    light rays bouncing off of a finite collection of    reflective wires (lines)  in $E = \R^3$.  
 
  \vskip .2 cm  
\noindent
\subsubsection{Motivating Example: $N$-body billiards}
\label{sss:motivating_eg}
   $E = (\R^d)^N$ is   the configuration space
 for the $N$ massive point particles moving in  $d$--dimensional Euclidean space $\R^d$. Endow $E$    with its mass metric, by which we mean
 the inner product whose squared norm is twice the kinetic energy.   
 Take  $\mL$  to consist  of   the   $N \choose 2$ binary collision subspaces 
\begin{equation}
\Delta_{ab} = \{ q = (q_1, \ldots , q_N), q_i \in \R^d: q_a = q_b\} \subset E
\quad \mbox{ for some pair }a \ne b.
\label{def:Delta:ij}
\end{equation}
We call this class of examples ``$N$-body billiards''. 
 See the next section for  details.

%{\coma [I am using $a, b$ instead of $i, j$ because of $q(t_i) = q_i$ . Okay?]}
 
\subsubsection{Billiard Rules.} 
\label{sss:rules}   We now define  what it means for a polygonal curve   $q: \R \to E$ to be a billiard trajectory. 
By a collision point for $q$ we mean a time $t$ or the corresponding point $q(t)$ such that $q(t) \in C$.  Thus at a collision point  $q(t) \in L$
for some $L \in \mL$.  
We assume that  collision points are discrete.  {\it In particular no edge of $q$  lies within an $L$.}  
Every vertex of $q$ is a collision point.  
The velocities  $v_-, v_+$ of $q$  immediately  before and after collision with $L \in \mL$ are well-defined and locally constant. 
They suffer a   jump $v_- \mapsto v_+$ at collision.
Let  
 $$\pi_L: E \to L  $$
 be   orthogonal projection onto $L$.   We require each  velocity  jump to obey the  rules:
 \begin{equation} 
 \|v_- \| = \| v_+\|  \qquad \text{ ``conservation of energy''}
 \label{eq:energy}
\end{equation}
  \begin{equation} \pi_L  (v_+) = \pi_L ( v_-) \qquad \text{ ``conservation of momentum''}
 \label{eq:momentum}
\end{equation}
 In N-body billiards (\ref{sss:motivating_eg})  these   rules correspond to conservation of energy and momentum. 
 (Note that the rules allow for no jump:  $v_- = v_+$.) 

To summarize,    a billiard trajectory for $(E, \langle \cdot, \cdot \rangle, {\mathcal L})$ is an oriented polygonal curve in  $E$ 
with   vertices on collision subspaces, and no edge of which lies within a collision subspace.
At each collision the velocity jump $v_- \to v_+$  obeys  the two rules (\ref{eq:energy}, \ref{eq:momentum}) above. {\it  Without loss of generality we will assume the curve's  speed is } $1$.  

\subsubsection{Multiple collisions.}  
\label{sss:multiple_collisions} The attentive reader will have noticed that  the law of reflection (\ref{eq:momentum}) is ambiguous if the collision point $q_*$ belongs to more than one $L$.
This ambiguity  is analogous to the problem of  trying to define standard billiard dynamics at the corner pocket of a polygonal billiard table in the plane. 
 We get around this ambiguity by agreeing to {\it choose}  one of the collision subspaces to which $q_*$ belongs and then using only  that subspace in implementing law (\ref{eq:momentum}).
 Thus we view   billiard trajectories with multiple collisions as coming with the extra structure of a   labelling of  collision points, with each collision point being 
 labelled by one of the   $L \in {\mathcal L}$ to which it belongs.  For more on problems arising with multiple collisions see subsection \ref{ss:McGehee} further on.

\subsubsection{Dimension and   Transversality.}   For simplicity of exposition {\it we will henceforth assume that
each subspace $L$ has the same codimension } $d$ and $d \ge 1$.    This 
assumption excludes various pathologies such as $L_1 \subset L_2$ occuring within our collection $\mL$  of subspaces.

In addition to being all of   the same codimension $d$, the collection  (eq (\ref{def:Delta:ij})) of collision subspaces for N-body billiards 
are  pairwise transversal: ${\rm codim}(L_1 \cap L_2) = 2d$
 for all distinct pairs $L_1, L_2 \in \mL$.  We believe such   transversality assumptions may be very useful in future work.

\vskip .2 cm  
\noindent
\subsubsection{Non-deterministic Dynamics.}  
\label{sss:nondeterministic}   
For a given incoming $v_- \in E \setminus 0$  to a $q_* \in L$ there is a    $(d-1)$--dimensional sphere's worth of choices for the outgoing $v_+$'s,
namely the set of all solutions to eq (\ref{eq:energy}, \ref{eq:momentum}) for that fixed $v_-$. 
It follows that the billiard process is  non-deterministic: there is no univalued rule that takes us from the past motion to the future motion.
However, we do not view our   billiard dynamics as a stochastic process.   Rather we think of 
our  billiard trajectories   as arising as   limits
of deterministic $N$-body dynamics, and we are interested in what is the set of all    possible limits.
%We simply know what it means to be a ``solution'' -that is a billiard trajectory, without specifying the `process'' by which these billiard trajectories arise.
See section  \ref{s:motivation} below.  

(Even if  $d =1$ we do not have deterministic dynamics, since the $0$-sphere consists of two choices.  
It is standard to turn this case into a deterministic dynamics by requiring {\it transversality}:  $v_+ \ne v_-$  at each collision.
This is what is done for $N$ point particles moving on the line:   the  dynamics preserves their order on the line. The game  is equivalent 
to playing billiards on a closed   polyhedral cone in $\R^N$.)

\subsection{Basic Questions and Main result.} 
\label{ss:MainResult}
\subsubsection{Fundamental Finiteness Theorem.}\quad\\ 
\label{sss:FFT}

QUESTION 1.  Is the total number of collisions of a billiard trajectory finite? 

 The answer is the fundamental theorem of the subject.
 
 \begin{backtheorem} [\cite{BFK1, BFK2}] 
 \label{thm:background}There is a $K= K( E, \langle \cdot, \cdot \rangle, {\mathcal L})$ such that every trajectory has less than or equal to $K$ collisions.
 \end{backtheorem} 
   To appreciate the subtlety of the problem of computing the smallest $K$,  even in apparently simple deterministic ($d=1$)  situations,  we strongly urge the reader to take a peek at \cite{CountingPi}.

\subsubsection{Itineraries.}  
 
\begin{definition}The {\em itinerary} of a billiard trajectory is the list of collision subspaces $L \in {\mathcal L}$ that  it intersects, in their   order of occurrence. 
\end{definition}

By the  Background  Theorem  of BFK  just stated,   any realized  itinerary has length less than or equal to $K$.  So
if there are a total of  $M$ subspaces in $\mL$, then the set of all realized itineraries is a finite set of length less than $M^K$.
(Repeats such as $L_1 L_1 \ldots$ are not allowed, hence the strict inequality.)   

QUESTION 2.  What is the finite set of all itineraries which are realized by some billiard trajectory?

This is a hard question about which we have very little to  say. 

We can observe that beyond $L_{i+1} \ne L_i$ there may be other
`topological' restrictions on the allowable itineraries.  For example, for 
 $N \ge 4$ bodies on the line ($d=1$; $E = \R^N$)  after the itinerary $(\Delta_{12},\Delta_{34},\Delta_{23})$
particles 1 and 4 can no longer  be neighbors, whether or not   collisions  change the ordering (are transverse). 
So  $(\Delta_{12},\Delta_{34},\Delta_{23},\Delta_{14})$ cannot be realized.

In  section \ref{s:examples} we give some partial results regarding this question when  each $L$ is a line.

\subsubsection{Space of trajectories realizing a given  itinerary}

Suppose   a particular  itinerary is realized.  We can then ask about all of its realizations.

QUESTION 3.  What is the structure (dimension, smoothness, symplectic character) of the space  of all    billiard trajectories having a given
itinerary?

%, nor can one contain the other.  

\vskip .2 cm  
\noindent
   
% with  universal upper  bound $C$ depending only on the Euclidean geometry of the arrangement  ${\mathcal L}$.  
% (The bound of \cite{BFK1} is almost certainly sub-optimal and is   superexponetial in the cardinality of ${\mathcal L}$.) 

Answering  Question 3 is the point of our paper.  
\vskip .15cm 
{\it From  now on we  fix an   itinerary}   $L_1 L_2 \ldots L_k$, with $L_i \in \mL$.  

{\it Defining the space $\mB$ of point billiard trajectories realizing the  itinerary.} 
Write  ${\mathcal B} (L_1 \ldots L_k)$, or simply ${\mathcal B}$
for the space of all billiard trajectories realizing the given itinerary.  Let us be more precise:
a  billiard   trajectory $q$ is in the subset $\mB (L_1 \ldots L_k)$ if and only if there are exactly   $k$ distinct collision times: 
\begin{equation}t_1 < t_2 < \ldots t_k,   \qquad q(t_i):= q_i  \in  L_i .
\label{intersections}
\end{equation} 
We  emphasize that 
\begin{equation}
\label{noncollision}
q_i \ne q_{i+1},\qquad i =1,2, \ldots, k-1 
\end{equation}
since $|q_{i+1} - q_i | = t_{i+1} -t_i  >0$.  
Also 
$q_{i+1} -q_i, q_i - q_{i-1} \notin L_i$ since no edge of $q$ lies in the collision locus.
Condition (\ref{noncollision}) does not exclude the possibility
of $q_i \in L'$ for some $L' \in \mL$, $L' \ne L_i, L_{i+1}, L_{i-1}$.  In this case
we   label $q_i$ with   $L = L_i$ when applying our   `conservation of momentum'' rule eq (\ref{eq:momentum}).
 We endow ${\mathcal B}$
with the compact-open topology.  

A trajectory  $q \in  {\mathcal B} (L_1 \ldots L_k)$ has an initial ray $r_-$    parameterized by the initial segment $(-\infty, t_1)$ 
where $q_1 = q(t_1)$ is the first collision.  Similarly $q$ has a final ray $r_+$    parameterized by the final segment $(t_k, \infty)$ 
where $q_k = q(t_k)$ is the final collision along $q$. 
Extend the rays to oriented lines  $\ell_-, \ell_+$.  We want to think of the fixing of the itinerary
as defining a ``scattering map''  
\begin{equation}
\ell_+ \mapsto \ell_+
\label{scattering}
\end{equation} on the space
$\LINES (E)$ of oriented lines in $E$. 
 (We will elucidate the structure of $\LINES(E)$ as a symplectic manifold momentarily.) 
However, this ``scattering map''  is almost never a map in that one $\ell_-$ may give
rise to many $\ell_+$'s.  See example \ref{ex:1} below.
Instead we have ``scattering relation''  
$$\mR = \mR (L_1 \ldots L_k)  \subset \LINES (E) \times \LINES (E).$$
\begin{definition}  The scattering relation $\mR = \mR(L_1 L_2 \ldots L_k)$ associated to  the chosen itinerary  $L_1 L_2 \ldots L_k$
 consists of all pairs   $(\ell_-, \ell_+)$ 
of incoming and outgoing lines for billiard trajectories $q \in {\mathcal B} (L_1 \ldots L_k)$.
\end{definition} 
We have just defined a  continuous map
$$\mB \to \mR \subset \LINES (E) \times \LINES (E)$$
which sends each  trajectory $q \in \mB(L_1 \ldots L_k)$   to its incoming and outgoing (oriented) lines.  The image of
this map is the scattering relation.  
The group of time translations acts on the space of billiard trajectories,  sending $q(t)$ to $q(t-t_0)$, for $t_0 \in \R$, 
without altering the itinerary or  the incoming or outgoing
line.   Thus our map into the scattering relation    induces  a map on the quotient domain with  the same image.  We name this map the {\it scattering projection}.
\begin{equation}
{\rm\it SCAT}: \mB/\R  \longrightarrow \mR \subset \LINES(E) \times \LINES(E)
\label{SCAT}
\end{equation}

We can now state our main result.
\begin{theorem} 
\label{main} 
The scattering relation $\mR$  is a  Lagrangian relation on the symplectic manifold  $\LINES (E)$
of oriented lines in $E$. In particular $\mR$ is a smooth manifold of dimension $2(\dim(E) -1)$. 
The scattering projection (eq (\ref{SCAT})) defines a diffeomorphism
between  ${\mathcal B}/\R$  and  $\mR$.
In   particular, modulo time translation,
a  point billiard trajectory  realizing the given itinerary
is uniquely  determined by its incoming and outgoing lines.  
\end{theorem}
For completeness, we recall for the reader the definition of ``Lagrangian relation'' and the symplectic structure on $\LINES (E)$ in what immediately
follows.

\subsubsection{Lagrangian relations} 
\begin{definition}  
A {\em Lagrangian relation} on a symplectic manifold 
 $(P, \omega)$ is a Lagrangian submanifold $\mR$ of the  product symplectic  manifold 
 $\bar P \times P$, where the bar of  ``$\bar P$'' means we endow the product with  the symplectic structure $-\omega \oplus \omega$. 
\end{definition} 
Graphs of symplectic maps $P \to P$ are Lagrangian relations.  We
think of Lagrangian relations as  generalized symplectic maps,
that is, symplectic maps  which are ``allowed to go vertical'' at various places.

\subsubsection{The symplectic structure on the space of lines.}
\label{sss:lines} 
An oriented line $\ell \in \LINES(E)$   can be 
represented by  an initial position  $A \in E$
and an initial velocity $v_A \in  E_v \cong E$. 
(We  use  the subscript  $v$ in ``$E_v$''   to keep track of who is a  velocity and who is a position.) 
The line associated to $(A, v_A)$ is     parameterized as $A + t v_A$.  
We will  insist that   velocities $v_A$ are unit: $|v_A |=1$.
$(A, v_A)$ and $(C, v_C)$ represent the same oriented line if and only if $v_A = v_C$ and
$C = A + s v_A$ for some real number $s$. 
There is a unique point $Q \in \ell$ closest to the origin of $E$. This $Q$ is 
determined by the algebraic condition   $\langle Q, v_A \rangle = 0$.  
Choosing   $Q$ as the initial position $A$ on $\ell$  sets up a diffeomorphism between the space $\LINES (E)$
of oriented lines in $E$ and 
 the tangent bundle of the unit sphere in $E$:
\[ \LINES(E)  \cong T S(E_v) =\{(v,Q) \in E_v \times E:  |v| =1,  Q \perp v \} . \] 
Use  the Euclidean structure to identify $TS(E_v)$ with $T^* S(E_v)$, 
thereby giving the space of lines a symplectic structure.

{\sc Remark.} 
The  diffeomorphism $\LINES(E) \to TS(E_v)$     reverses  the role of positions 
and velocities.  The position 
$v \in S(E_v)$ at which the tangent vector $(v,Q)$ is attached  represents the velocity  
vector $v = v_A$ of the corresponding line,
while the tangent or $Q$-part of $(v,Q)$ represents an initial  position  point on the line $\ell$, namely the closest point to $0$.

\subsubsection{ Lines as a reduced space} 
\label{sss:reduced}

  The space of oriented lines  can be recast as   a  symplectic reduced space. 
    Let $H (A, v) = \frac{1}{2} |v|^2$ be the usual
Hamiltonian for free particle motion. Here   $(A,v) \in E \times E_v \cong E \times E^*  \cong  T^*E$.    
 The flow of  the Hamiltonian vector field for $H$ is  $\phi_t (A, v_A) = (A + tv_A, v_A)$ which  is a symplectic  $\R$
action on the full  phase space.  Its integral curves are  lines.
 The   level set $H^{-1}(1/2)$ consists of 
 those   initial conditions $(A,v)$ such that $|v| =1$.  The space $\LINES (E)$  of oriented lines is thus the  sub-quotient $H^{-1}(1/2)/\R$ of  $E \times E_v$  
by this $\R$ action.  This sub-quotient construction is  precisely the symplectic reduction construction: $\LINES(E)$
 with its symplectic structure  is an instance of the construction of 
the ``symplectic reduced space''.     Write
\begin{equation}
\pi:  E \times S(E_v) \to \LINES (E) 
\label{reduction}
\end{equation}
for the corresponding quotient map.  Thus
$\pi (A,v) = \pi (\tilde A, \tilde v)$ if and only if $\tilde v = v$ and $\tilde A = A + tv$ for some $t \in \R$. 

\subsubsection{ The unreduced scattering relation} 

In order to  prove and to  better  understand   our main   theorem \ref{main} we must    ``unreduce'' the relation $\mR$ by  working  directly with  normalized 
initial  conditions $(A, v_A) \in E \times S(E_v)$ instead of the associated   oriented line $\pi(A, v) = \ell$.  
If $q(t)$ is a billiard trajectory in ${\mathcal B} (L_1 \ldots L_k)$ consider again its initial ray $r_- \subset \ell_-$ and final ray  $r_+ \subset \ell_+ $.
Pick corresponding points $A \in r_-$,  $B \in r_+$ and the corresponding directions $v_A, v_B$.
{\it We emphasize that we are saying nothing about the times $t_A, t_B$ at which the   points $A, B$ are selected
along $q(t)$.}
In this way we have chosen pairs $(A,v_A), (B, v_B) \in E \times E_v$.   
The unreduced  statement of  theorem \ref{main} is 

%\eject

\begin{theorem}  \label{main_unreduced} For each  $q \in \mB$ consider the two-parameter family  of  pairs of 
boundary conditions   
$$\big((A,v_A),  (B, v_B)\big) = \big((q(t_0), \dot q (t_0)),  (q(t_{k+1}), \dot q (t_{k+1}))\big)$$
lying on  the incoming and outgoing rays of $q$.  (Here $t_0 <t_1$ and $t_{k+1} > t_k$
as per eq (\ref{intersections}).) 
As $q$ varies over $\mB$ these pairs 
%$((A,v_A),  (B, v_B))$ 
sweep out a Lagrangian relation $\tilde \mR = \tilde \mR (L_1 \ldots L_k)$
on $E \times E_v$. The projection $((A, v_A), (B, v_B)) \mapsto (A, B)$ maps 
 $\tilde \mR$ diffeomorphically onto  an open subset
of $E \times E$.  The projection  of $\tilde \mR$ to $\LINES(E) \times \LINES (E)$
by $\pi \times \pi$ (where $\pi$ is as in eq \ref{reduction})  is the relation $\mR$ of theorem \ref{main}. 
\end{theorem} 
 
\vskip .2 cm  
\noindent
{\sc Remark on algebraicity.}  
Our Lagrangian  relations are semi-algebraic varieties: they are defined by 
algebraic equations together with algebraic inequalities.  This fact  follows from our proof of the theorem using generating functions.

{\sc Remark: Scaling, Symmetries and Conservation Laws. }  
Point billiard trajectories  enjoy a scaling symmetry. $N$-body billiards
 enjoy translational and rotational symmetries and the consequent conserved quantities
of linear and angular momentum.  Details  of these symmetries
are discussed in section \ref{s:symmetries}.  
\noindent

\section{Motivation : The Gravitational $N$--Body Problem.}   
\label{s:motivation}

We go into some detail regarding our underlying motivation.
The basic set-up, $E = (\R^d)^N$ with the collision subspaces
being the binary collision subspaces $\Delta_{ab}$ was described above in   subsection \ref{sss:motivating_eg}
and we keep the same notation.

Positive energy solutions to the  gravitational two-body problem,     
 viewed in a center-of-mass frame, consist of a pair of  coplanar hyperbolas sharing the   origin as a  focus.  Viewed from   afar away,  these hyperbolas become indistinguishable from their  asymptotes:
 the two bodies come in along their separate rays, bounce off each other, to head back to infinity along different rays.     
 %By    homogeneity of the kinetic and the potential energy, moving far away is the same as making the energy large.    
 
For  the gravitational $N$-body problem the same space-time  picture holds  when viewed   from away from all close encounters.  
Each body moves   nearly on  a straight  line at nearly  constant speed until it   comes  into very close vicinity of another body   at which time  it veers off 
to recede along another near-line at near-constant speed.  In the limit \footnote{ The limit is $\frac{q(\lambda t)}{\lambda}$ as $\lambda \to \infty$ where $q(t) \in (\R^d) ^N$}, what happens at  these close encounters is 
the bodies ``bounce off'' each other.  The direction of this  ``bouncing'' will look   random unless we know   detailed specifics
of the incoming motion.   Without these details, all we can say  is that each  bounce is an elastic collision :    total energy and linear momentum
are conserved.  These two conservation laws are   encoded by our rules of reflection (eq (\ref{eq:energy}, \ref{eq:momentum}).

Thus we expect
%\footnote{ The limit is $\frac{q(\lambda t)}{\lambda}$ as $\lambda \to \infty$ where $q(t) \in (\R^d) ^N$} 
certain families of   positive energy solutions to the graviational  $N$-body problem will limit onto    $N$-body billiard trajectories as 
 described above (see subsection \ref{sss:motivating_eg}).
 % with the collision subspaces
%being the binary collision subspaces denoted  $\Delta_{ab}$  above (see eq \ref{def:Delta:ij}). 
%$\R^d$ is the underlying $d$-dimensional Euclidean space   which the $N$ bodies move.
%The collection of linear subspaces ${\mathcal L}$ are the $N \choose 2$ binary collision subspaces  $\Delta_{ab}$ as given in eq \ref{def:Delta:ij} above. 
%The Euclidean structure on $E$ is the ``mass metric''  $| dq | ^2 = \Sigma m_a |dq_a|^2$ where  $m_1 , \ldots, m_N$  are the masses and $|dq_a|^2$ is the standard Euclidean
% structure on $\R^d$.     The condition of eq (\ref{eq:momentum})   is precisely  conservation of linear momentum of the total system
% at   collision.  Conservation of total energy is encoded in the requirement that the trajectories have unit speed. 
  In a subsequent paper we will   prove this assertion by showing that   $N$-body billiard trajectories  are  ``shadowed'' by  families of trajectories
   of positive energy solutions to  the gravitational $N$-body problem. 
 
 \vskip .2 cm  
\noindent
\subsection{Multiple Collisions and clusters.} 
\label{ss:McGehee}
A collision between   three or more particles (or two or more simultaneous binary collisions) 
corresponds to a point $q \in (\R^d)^N$ lying in several $\Delta_{ab}$.  
The paper of Mather and McGehee \cite{McG}, and subsequent work on non-collision singularities based on their ideas 
make suspect the validity of our  underlying assumption 
(\ref{eq:energy}) of conservation of kinetic energy when trying to model such   multiple  collision events with point billiards. 
Mather and McGehee establish the existence
of a set of initial conditions for 4 bodies (on the line) where the kinetic energy starts out $O(1)$
and in a finite time becomes arbitrarily large, arbitrarily far away from the close encounter region.
   The infinite negative potential energy well of near-triple collision serves as a source
which one of the bodies can extract to make its speed arbitrarily high.   We imagine 
 the following caricature of celestial mechanics based on the notion of cluster decompositions \cite{DG},
 where each  clusters represents a  subset  of close   tightly bound particles.   
Total energy and momentum is preserved for each isolated cluster.   
But not all energy need be   kinetic. We could even allow trajectories to  move inside intersections of the $\Delta_{ab}$, corresponding
to systems that are bound over some large interval of time. At collisions
between clusters, corresponding groups of particles can experience inelastic
scattering,   potential energy being stored in  groups or released from it, and   redistributed.
\vskip .2 cm

\section{Generating Families and the proof.}
\label{s:generating_families}

 The chord length between
successive impacts of the   ball with the table serves as the generating function for the
standard billiard map associated to a convex table in the plane.  So it is
not a great surprise that the   path length of  finite segments of polygonal paths realizing the given itinerary  serves a similar   function
for our non-deterministic billiard processes. 
% A minimization problem underlies our proof of the theorem and our   perspective  on the problem. 
Fix points $A= q_0$ on the incoming ray and   $B=q_{k+1}$
on the outgoing ray of the  billiard trajectory $q \in \mB$.
Let $q_i \in L_i$ be the intermediate collision points as per eq (\ref{intersections}).
Then the length of the segment $q([t_0, t_{k+1}])$ is: 
\begin{equation}
S(A, q_1, \ldots, q_k, B) = |A-q_1| + |q_1 -q_2| + \ldots + |q_{k-1} - q_k| + |q_k -B|
\label{action}
\end{equation}
and this is  also the travel time of this segment. We turn this observation around
to find the billiard trajectories as critical points of $S$.   

\vskip .2 cm  
\noindent
{\sc Minimization Problem.} \\  
Fix $A, B \in E \setminus C$.  Minimize (\ref{action})  over all
intermediate choices   $q_i \in  L_i$.
%over all $k$-tuples $(p_1, p_2, \ldots, p_k) \in L_{B_1} \oplus L_{B_2} \ldots L_{B_k}$.

%The global minimum depends continuously on $A,B$. Let us call a minimum ``transverse''
% if $p_i \ne p_{i+1},  I =0, \ldots, k$.  ..

\vskip .2 cm  
\noindent

 Write $xy$ for the line segment joining $x$ to $y$, $x, y\in E$, parameterizing $xy$
so as to have unit speed.   If $x_i \in E$ are a collection of points then by $x_1 x_2 x_3 \ldots x_n$ we will mean the polygonal path
with $n-1$  edges $x_i x_{i+1}$.
 Let
$$ \Lambda = L_1 \times L_2 \times \ldots \times L_k.$$
Then if $A, B \in E$  and $\lambda = (q_1, \ldots , q_k ) \in \Lambda$
we  write $A \lambda B$ for the piecewise linear segment $A q_1 q_2 \ldots q_k B$.

\begin{definition}
\label{def:generic} We  will say that $\lambda \in \Lambda$ is ``generic''
if  $q_1 \notin L_2, q_{k} \notin L_{k-1} $ and $q_i \notin L_{i-1} \cup L_{i+1}, 1 < i < k$. 

\noindent We will say that  $(A, \lambda, B) \in E \times \Lambda \times E$ is `generic'' if $\lambda$ is generic,  if $A, B \in E^0$
and if the rays $q_i A$ and $q_k B$ have no collisions besides their initial points $q_1, q_k$.
%(Recall from eq \ref{table} that $E^0$ is our ``table'' , which simply means that  neither $A$ nor $B$ are collision points.)

\noindent We will call the open set of all generic points in $E \times \Lambda \times E$ the ``generic set''. 
\end{definition} 

Define  
\begin{equation} S_{A,B} : \Lambda \to \R;\quad 
S_{A,B} (q_1, q_2, \ldots, q_k) = S(A, q_1, q_2, \ldots, q_k, B)
\label{actionAB}
\end{equation}
by viewing the action (\ref{action}) to be a function of the intermediate
intersection points  $q_i$ alone, with $A, B \in E$ as parameters.

\begin{proposition} \label{critical_pts}
%\begin{itemize} \item{1)} 
Suppose   $(A, \lambda, B)$ is generic in the sense of definition \ref{def:generic}.  Then the following are equivalent.\\
\begin{itemize}
\item{}(A)  $\lambda $ is a critical point of $S_{A,B}$\\
\item{}(B) $A \lambda B$ is a segment of a  billiard trajectory realizing the given itinerary.
\end{itemize} 
If either  condition holds  then the direction of the incoming line of the associated billiard trajectory   is 
$v_A = - \nabla_A S (A, \lambda, B)$ while the direction of the outgoing line is 
$v_B = + \nabla_B S (A, \lambda, B)$ where 
$\nabla_A S$, $\nabla_B S:  E^0\times \Lambda \times E^0 \to E$ 
denote the gradients with respect to the $A, B$ variables.  
\end{proposition}

\begin{example} \label{ex:1} [Total Collision]
Consider  the case $\mL = \{ 0 \}$, so that  the only subspace  is the 0 subspace. 
A  linear billiard trajectory   realizing the itinerary $(0)$ consists of an angle with vertex at $0$. 
   The parameter space $\Lambda$ is the single point $0$.  The action  is
$S(A, 0, B) = |A| + |B|$.
The intermediate collision point   $q_1 = 0$   cannot be varied  so   the condition $d_{\lambda} S = 0$ is  vacuous.
We compute   $\nabla_A S  = A/|A|, \nabla_B S = B/|B|$
consequently the Lagrangian relation of theorem \ref{main_unreduced}  
consists of all quadruples $((A, v_A),( B, v_B)) \in  TE \times TE$
for which $v_A =   -A/|A|$ and  $v_B = B/|B|$, and $A, B \ne 0$. 
The first pair $(A, -A/|A|)$ represents  the initial position and velocity of a line thru the origin,  
moving towards the origin.
The final pair $(B, B/|B|)$ represents  the initial position and velocity for  a line thru the origin moving away from  the origin.
Our incoming line and outgoing  line both pass through the origin, so their  
``$Q$ parts'' are   $0$.  (See subsubsection \ref{sss:lines}.)   Their $v$ parts, $v_A$ and $v_B$ are arbitrary unit vectors. 
The Lagrangian relation $\mR$ of theorem  \ref{main}  is 
 the product of the two zero sections of $T^* S(E) = \LINES(E)$.
\end{example}

Proposition \ref{critical_pts} asserts that  $S$ is a  ``generating family''  (also known as a ``Morse family'') for the   Lagrangian  relation of theorem \ref{main_unreduced}.
We recall  the notion of a  generating family.  
\begin{definition} 
\label{def:generating_family}
The function  $F: E \times \Lambda \times E \to \R$ is a {\bf generating 
family}  for  the Lagrangian relation $\mR$ on $E \times E_v$  if $\mR$ consists of   those quadruples  (pairs of pairs) 
 $\big((A, v_A),(  B, v_B)\big) \in (E \times E_v) \times (E \times E_v)$ for which  there exists a 
 $\lambda \in \Lambda$ such that 
 \begin{itemize}
\item{} (i) $(A,\lambda, B)$ is a smooth point of $F$, and 
\item{} (ii) 
$d_{\lambda} F (A, \lambda, B) = 0$,  $v_A = -\nabla_A F (A, \lambda, B)$
and $v_B = + \nabla_B F (A, \lambda, B)$.
\end{itemize} 
Here $\nabla_A, \nabla_B$
are the gradients with respect to these first and last component variables, $A, B$ and  
$d_{\lambda} F(A, \lambda, B) \in  \Lambda^*$ is the differential with respect to $\lambda$. 
\end{definition}  
\vskip .2 cm  
\noindent
%\subsubsection{Historical Remark.} 

The notion of generating family was  formalized
by H\"ormander in \cite{Hor1} \cite[Def.\ 25.4.3]{Hor2} under the name of  ``phase function''.   Libermann and Marle 
\cite{LiMa}   use the name  ``Morse family'' and we find their treatment exceptionally clear. 
({\sc See  Definition 1.10 in \cite[Appendix 7.1]{LiMa}.}) Paraphrasing: 
``Let $\pi : B \to N$ be a submersion and $S: B\to {\mathbb R}$ be a 
differentiable function. 
The function $S$ is called a {\em Morse family} (for $N$, or for $R \subset T^*N$)  if  the image of the one-form
$dS: B \to T^*B$  and the conormal bundle to the fibers of $\pi$, are transverse within $T^*B$.  This transverse intersection
is necessarily smooth and  pushes down to $T^*N$
where it forms a  Lagrangian submanifold $R$, the Lagrangian submanifold for which  $S$ is the `Morse family'.''   

The transversality condition in the definition just given  of a Morse family is needed to insure that the corresponding
Lagrangian submanifold  is smooth.  In our case we establish 
smoothness by establishing:

\begin{proposition} \label{Hessian}
Every  critical point $\lambda$ of $S_{AB}$  which is a generic point  in the sense of definition \ref{def:generic} is a non-degenerate
critical point, so transversality holds as discussed above.  Indeed, the   Hessian of $S_{A,B}$ at $\lambda$ is positive definite.
% If $\lambda = (p_1, p_2, \ldots, p_k)$ is a critical point of $S_{A,B}$  which is transverse relative to $A, B$   then $A, \lambda, B$ is a smooth point of $S$ and the
\end{proposition}

%{\color {red} bit of reference, history, on generating families here??}

\subsection{\sc Proof of proposition \ref{critical_pts}. } 
\label{ss:critical_pts}
 
 For the function $x \mapsto |x|$ we have that $d|x| = \frac{ \langle x, dx \rangle} {|x|}$.
 (The  algebraic meaning of `$dx$' here, as per computations found frequently in  Chern or Cartan, is that $dx$ is the  identity map on $E$, this being the differential of the map $x \mapsto x$.
In other words, for $v \in E$   $(d|x|) (v) = \frac{ \langle x, v \rangle} {|x|}$.) Similarly if $x_0 \in E$ is a constant vector then
 $d |x - x_0| = \frac{ \langle x-x_0, dx \rangle} {|x-x_0|} = \langle n(x, x_0) , dx \rangle$, where we
 write $$n(x,y) = \frac{ x-y}{|x-y|}$$
  for the unit vector pointing from $y$ to $x$, assuming $x \ne y$. 
   Now write $d_i$ for the differential of $S_{A,B}$ with respect to $q_i$, keeping the other $q_j$ constant.
 We have
 $$d_i S_{A, B} = d_i (|q_{i-1} - q_i| + |q_i - q_{i+1}|) =  \langle n(q_i, q_{i-1}), dq_i \rangle +  \langle n(q_i, q_{i+1}), dq_i \rangle . $$
 Since $n(y,x) = -n(x,y)$, this yields
 $$d_i S_{A, B} = \big\langle n(q_i, q_{i-1})- n(q_{i+1}, q_i)\ ,\ dq_i \big\rangle . $$
 Now $dq_i$ is the identity on  $L_i$, so this differential is zero if and only if
 $n(q_i, q_{i-1})- n(q_{i+1}, q_i) \perp L_i$ which is the same as requiring that 
 $\pi_i (n(q_i, q_{i-1})- n(q_{i+1}, q_i)) =0 $, where we have written $\pi_i$ for $\pi_{L_i}$.
But if the piecewise linear trajectory $A  q_1 q_2 \ldots q_k B$
is parametrized by arc length, traveling from $A$ to $B$, then its
velocity just before collision with $L_i$ is $v_{i, -} = n(q_i, q_{i-1})$
and its velocity just after collision is $v_{i, +} = n(q_{i+1}, q_i)$, so that our condition
of criticality is equivalent to the condition of conservation of momentum (equation (\ref{eq:momentum}))  at collision $i$.
Finally $dS_{A,B} = 0$ if and only if for $i = 1, 2, \ldots, k$ we have
$d_i S_{A, B} = 0$.  \hfill $\Box$\\
\vskip .2 cm  
\noindent

We postpone  the proof of proposition  \ref{Hessian} to section \ref{s:Hessian}.

%{\color {red}: maybe turn this into proof of theorem 1?} 

\vskip .2 cm  
\noindent
\subsection{Proof of (most of)  theorem \ref{main_unreduced}.} 

%Suppose that ${\mathcal B} (L_1 \ldots L_k)$ is non-empty. 
Let $q_0\in {\mathcal B} (L_1 \ldots L_k)$ with initial ray $\ell_{in, 0}$ and final ray $\ell_{out, 0}$.
 Let  $\lambda_0 = q_1 ^0,  q_2 ^0 ,\ldots,  q_k ^0$ be its
collision points.  According to the definition of a billiard trajectory we cannot have $q_{i +1} \in L_i$, $1 \le i \le k-1$
for otherwise  segment $q_i q_{i+1} \subset L_i$ which is forbidden.   Similarly $q_{i-1} \notin L_i$ for $2 \le i \le k$.
Choose points $A_0  \in \ell_{in, 0}$,    $B_0 \in \ell_{out, 0}$. Then  $A_0, B_0 \notin C$.
Thus $(A_0, \lambda_0, B_0) \in E \times \Lambda \times E$ is a  generic  point.  
And according to proposition \ref{critical_pts},  $\lambda_0$ is a  critical point 
of $S_{A_0,B_0}$.

Now flip the logic around. 
% where $A_0$ is a point on the initial ray of $q_0$ and $B_0$ is a point on the final ray. 
%, as is every other point of ${\mathcal B} (L_1 \ldots L_k)$.
Consider the map $E \times \Lambda \times E \to \Lambda^*$ 
\begin{equation}
\label{crit_pt_map}
(A, \lambda, B)\  \longmapsto \ dS_{A,B} (\lambda) \in \Lambda^* .
\end{equation}
Proposition \ref{critical_pts}  asserts that  the zeros $(A, \lambda , B)$ of this map  which are   generic points
(in the sense of definition \ref{def:generic}) are  precisely   the  billiard segments for some $q \in \mB (L_1, L_2, \ldots, L_k)$. 

The  chosen segment of $q_0$ from $A_0$ to $B_0$ is such a zero.   
We   use Proposition \ref{Hessian} in conjunction with the  Implicit Function Theorem to get nearby, smoothly varying, zeros. 
The derivative of  map (\ref{crit_pt_map}) with respect to $\lambda \in \Lambda$ at $A_0 \lambda_0 B_0$ is the Hessian of $S_{A_0, B_0}$ with respect to $\lambda$, evaluated at $\lambda_0$.
Proposition  \ref{Hessian} asserts this derivative is invertible.  The hypotheses of the Implicit Function Theorem hold.  There exist neighborhoods $U_- \subset E$ of $A$
and $U_+ \subset E$ of $B_0$ and 
a smooth function $U_- \times U_+ \to \Lambda$, 
written $(A, B) \mapsto \lambda(A, B)$    such that 
$A \lambda (A,B) B$ is a zero of the map (\ref{crit_pt_map})  and hence potentially part of a  billiard segment lying in $\mB$.
We complete this  billiard  segment to a full trajectory $q: \R \to E$ by  extending its initial and final segments $Aq_1$ and $q_k B$ to rays.  We can guarantee
that this extended full trajectory has no new collisions by taking a  sufficiently small neighborhood $U_-, U_+$ of $A_0, B_0$ and recalling that the generic
set is open.  This   full trajectory is now a billiard trajectory $q \in \mB$ with these $q$'s smoothly parameterized by their `endpoints''  by $(A,B) \in U_- \times  U_+$.

We   have  just described    billiards $q \in \mB$ as {\it locally} forming  graphs over their `endpoints' $(A,B)$.
By direct computation the velocity of the initial ray at $A$ is $v_A = (q_1 - A)/|q_1 -A| = - \nabla_A S$
while the velocity of final ray is $v_B = (B-q_k)/ |B - q_k| = \nabla _B S$.  Hence, when viewed in terms of initial and
final conditions  $((A, v_A),(B, v_B))$ at points along initial and final rays,  the space of  billiard trajectories $q \in \mB$  is realized {\it locally} as a Lagrangian relation $\tilde \mR$ on $E \times E$ which arises from   the generating family $S$, and forms
{\it locally}  a graph over some open set in $E^0 \times E^0$. 

It remains to prove that these local graphs piece together to a global graph over an open dense subset of
the space of endpoints  $E^0 \times E^0$.
That `piecing together' is precisely the uniqueness assertion of the penultimate sentence of theorem \ref{main_unreduced}
which states that a   billiard  trajectory $q \in \mB$ is uniquely determined (modulo
 time translations)  by its   endpoints $A,B$.  
Proving this uniqueness requires a new tool, summarized in  theorem \ref{min_thm} below. 

The assertion of the  last sentence  of   theorem \ref{main_unreduced} concerns the relation between the Lagrangian relation of theorem 2 and the relation described by  theorem 1.
The proof of this assertion is the same as  the proof of theorem \ref{main} which now follows.  
\hfill $\Box$

\subsection{Proof of  theorem \ref{main}: Reducing  Lagrangian Relations.} 
\label{sss:reducing_relations}
\vskip .2 cm  
\noindent

 We will  push the Lagrangian relation $\tilde \mR$ on $E \times E_v$ of theorem \ref{main_unreduced}   down to a Lagrangian
relation on   $\LINES(E)$ and verify that it is the desired Lagrangian relation $\mR$.

Recall from subsubsection \ref{sss:reduced} that $\LINES(E)$  is the symplectic reduced space
of $E \times E_v$ by the  Hamiltonian flow for the free particle Hamiltonian  $H(q, v) = \frac{1}{2} |v|^2$.  As such $\LINES(E)$  is a subquotient of $E \times E_v$
with subquotient map written $\pi: H^{-1} (1/2) \to \LINES(E)$.   
Observe that $\tilde \mR \subset  H^{-1} (1/2) \times H^{-1} (1/2)$,
since  whenever $((A,v_A), (B, v_B)) \in \tilde \mR$ then $v_A, v_B$ have unit length. 
Regardless of what points  $A$,  $B$ we pick  along the initial
ray $\ell_-$ and  final ray $\ell_+$  of a  fixed   billiard trajectory $q$, we get the same intermediate points $\lambda =  q_1 q_2 \ldots q_k$.
In other words,   $(A, v_A; B, v_B) \in \tilde \mR$ and $(A + h  v_A, v_A; B + s v_B, v_B) \in \tilde \mR$ give rise
to the same trajectory $q$, modulo translation (provided that  $h, s \in \R$ are appropriately restricted so we have
not ``passed'' the first or last collision of the initial or final ray).     In other words  these different choices
of $A, B$ yield the same initial and final rays, and hence the same initial and final lines.  
But this action of $(h, s) \in \R \times \R$ generates precisely  the kernel of  the   form $-\omega \oplus \omega$ on $(E \times E_v) \times (E \times E_v)$
upon restricting this  form to   $H^{-1} (1/2) \times H^{-1} (1/2) \subset (E \times E_v) \times (E \times E_v)$.    It follows that $\tilde \mR \subset  H^{-1} (1/2) \times H^{-1} (1/2)$ descends by  the  quotient  map
$\pi \times \pi$
to yield our desired Lagrangian relation $\mR \subset \LINES(E) \times \LINES(E)$. $\Box$.

%, yielding our desired Lagrangian relation $\mR = \mR (L_1 L_2 \ldots L_k)$.  

\vskip .2 cm  
\noindent

\subsection{Uniqueness. What remains to do. }   We    have established that our  space $\mB /\R$ of    billiard trajectories  realizing  the given itinerary, modulo time translation, is  {\it locally}  
 a graph over its initial and final rays.
But    theorem \ref{main_unreduced} and theorems \ref{main}  asserts that $\mB /\R$ is  globally a graph:   there cannot    be two billiard  trajectories with the given itinerary which 
share the same   initial and final rays.   The uniqueness assertion will
be established  by proving:  
\begin{theorem} 
\label{min_thm}  
For $A, B \in E^0$ 
there is a unique global minimum $\lambda \in \Lambda$  for $S_{A,B}$ and no other
 critical points or local minima. 
\end{theorem}

{\sc Caveat.} The global minimizer $\lambda \in \Lambda$ of  theorem \ref{min_thm}  might {\bf not}  yield 
a trajectory  in $\mB$
because it might    suffer multiple   collisions of the form  $q_i = q_{i+1}$
which were explicitly excluded from being paths in $\mB$.  See eq (\ref{noncollision}). 
Note that 
 $S_{AB}$ fails to be smooth at  such multiple collision  points. 

\noindent
The proof of theorem \ref{min_thm} will be given in section \ref{sect:CAT}.
\vskip .2 cm  

\noindent

{\sc Finishing the proof of the main theorem \ref{main}, given theorem \ref{min_thm}.} \\ 
Let $q \in \mB$ be a  billiard trajectory realizing the given itinerary. 
  Choose a point  $A$ on its initial ray,  $B$ on its final
ray, and let $\lambda \in \Lambda$ be the list of   collision  points  ticked off the itinerary.   By proposition  \ref{critical_pts},    
$\lambda$ is  a critical point for  $S_{A,B}$.   By  theorem~\ref{min_thm},   $\lambda$ is the global minimum of $S_{A,B}$ and  its only  critical point.
By proposition\ref{critical_pts} again,  there are no  other billiard trajectories  which pass
through $A$, tick off  the  given itinerary through a collision sequence, and then pass through $B$.  In particular  no other  billiard
trajectory  shares $q$'s   itinerary while having  the same initial and final ray.   
 This yields the uniqueness assertion of theorem \ref{main}  
 and the diffeomorphism assertion of the penultimate sentence of   theorem \ref{main_unreduced}. 
\hfill $\Box$

\section{The Gluing of CATs. Proof of theorem \ref{min_thm}.} \label{sect:CAT}

We follow  the non-smooth metric geometry  ideas and construction in \cite{BFK1} to obtain  the proof of theorem~\ref{min_thm}.
(See also \cite{BFK2}.)

\begin{figure}
\begin{center}
\includegraphics[height=4cm]{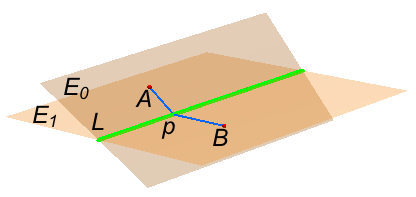}
\vspace{-5mm}
\end{center}
\caption{The metric space $E_L$}
\label{fig:Two:Planes}
\end{figure}

The main idea comes through clearly   by looking into the case of   a single subspace
$\mL  = \{L\}$, $L \subset E$, which is to say, an itinerary of length one.   
We form a new metric space $E_L$ by gluing two copies of
$E$ together along $L$.  We call the two copies ``sheets''  and label them 
$E_0, E_1$.  Thus
$$E_L = E_0 \cup_L E_1 ; \quad E_i \text{ copies of   } E.$$
See figure \ref{fig:Two:Planes} for the case where $L$ is a line in the plane.

We  define the metric on $E_L$ in terms of the minimizing  geodesics between two  points.
If the two points   lie in the same sheet   then the
geodesic between them is simply the usual  line segment  of  $E$ which joins them, 
viewed as 
lying in their shared  sheet.  If the two points $A$ and $B$ lie in different sheets,
the only way we can travel from $A$ to $B$ is by passing through   $L$ to cross from one sheet to another.  
We are led to the   problem of minimizing the distance from $A$ to $B$, in $E$,  among all
paths which touch $L$ in between.    In other words, we must minimize 
$S_{AB} (p) = |A -p| + |p-B|$ over $p \in L$.   As we have seen, there is a  unique minimizer $p_* \in L$.
Then the geodesic $AB$ consists of the union of the line segment  $Ap_*$ in $A$'s sheet
and the segment $p_* B$ in  $B$'s sheet.  The 
 minimization problem is the same problem we encountered earlier.  The geodesics
in this case are in bijection with the point billiard trajectories having  itinerary $(L)$. 

This construction of $E_L$  is a special case of a general  metric gluing procedure  which is the subject of a 
  theorem by Reshetnyak.   In order to describe Reshetnyak's theorem we must   recall   what it means for a metric space to
  be ``CAT(0)'' and ``Hadamard''.  
  %CAT(0) spaces are also called ``spaces with curvature bounded above by zero''.
    
  Let $(X,d)$ be a  path-connected metric space. 
    We  define the {\em length} of a path  in $X$ by taking infimums
    of ``polygonal approximations''  to the path (see \cite[Def. 2.3.1.]{Burago}).  
      $X$ is called a {\em length space}  if  
     the distance $d(A,B)$ between   points  $A, B \in X$ is  realized as the infimum of 
     the lengths
     of paths between $A$ and $B$.  If $X$ is also complete, then there is a  shortest such path, 
     denoted   $AB$,  and its length is $d(A,B)$.   
     We call $AB$ a {\em geodesic segment}.
  (There may be more than  one  geodesic segment joining $A$ and $B$.)  
  A {\em triangle}  $\Delta ABC$ in  $X$ is a subset consisting  of three points 
  $A, B, C \in X$ together with   geodesic segments  $AB, AC, BC$ joining them.
  A {\em Euclidean comparison triangle} $\Delta \bar A \bar B \bar C \subset \R^2$ for 
  $\Delta A B C$ is a triangle in
  the  Euclidean  plane whose sides are congruent to those of $\Delta ABC$:  
  $d(A,B) = \| \bar A - \bar B\|$,
  $d(A,C) = \| \bar A - \bar C \|$ and $d(A, C) = \| \bar B - \bar C \|$. 
If $x \in AB$  is a point on  side $AB$, then
there is a unique  comparison point $\bar x \in \bar A \bar B$ defined by    
$\| \bar A - \bar x \| = d(A,x)$ ,  $\| \bar x-  \bar B \| = d(x, B)$.

\begin{definition}  \cite[Defs. 4.1.9, 9.2.1.]{Burago})
A {\rm CAT(0) space} -- also known as a {\em space with non-positive curvature} -- is a complete length  space  $(X,d)$  
such that every sufficiently small  triangle $\Delta ABC$   in $X$  satisfies the following triangle comparison
property. Let $\Delta \bar A  \bar B  \bar C$ be a Euclidean comparison triangle for $\Delta ABC$.  
If
 $x \in AB$  
and  $\bar x \in \bar A \bar B $   is   the comparison  point then
 $d(x,C) \le \| \bar x - \bar C \|$. \\
\end{definition}

\begin{definition}  \cite[Defs. 4.1.9, 9.2.1.]{Burago}
A {\em Hadamard space} is a simply-connected CAT(0) space. 
\end{definition}

The CAT(0)  condition generalizes  the Riemannian geometry condition that  all sectional curvatures are non-positive
to the case of (possibly)  non-smooth metric length spaces.   
The fiducial example of a Hadamard space  is a Euclidean space.  Hyperbolic space and metric trees are other examples. 
\begin{theorem} {\bf (Reshetnyak \cite[9.1.21.]{Burago})}
If $X_1$ and $X_2$ are Hadamard spaces
containing isometric copies of the same convex set $K$,
then  the length space $X_1 \cup_k X_2$ constructed  by
gluing $X_2$ to $X_1$ along $K$ is again a  Hadamard 
space. 
\end{theorem}

%For a proof see \cite{Burago}, theorem 9.1.21. 

Reshetnyak's theorem asserts that the output of gluing 
Hadamard spaces  serves as another input!  We can   iterate.  
If $L_1 L_2 \ldots L_k$ is an itinerary we thus form the Hadamard ``itinerary'' space
$$E_{L_1 L_2 \ldots L_k} = E_0 \cup_{L_1}  E_1 \cup_{L_2} \cup \ldots \cup_{L_k} E_k.$$
Point billiard trajectories having  itinerary $L_1 L_2 \ldots L_k$   yield geodesics  which connect the first sheet $E_0$ to the last
sheet $E_k$.   

{\it Caveat.}  There may  be minimizing geodesics in this Hadamard itinerary space
which are not point billiard trajectories.  These will be  minimizers of $S_{AB}$ having either  a
multiple collision point $q_i = q_{i+1}$ or   an edge lying within an  $L_i$. 

\vskip .2 cm  
\noindent
{\sc Proof of theorem \ref{min_thm}.}  When   $X$  is   Hadamard
 there is a {\it unique} geodesic  $AB$ joining any two  points $A, B \in X$.
 %Moreover this   geodesic depends continuously on its endpoints.  
 See, for example,  Theorem 9.2.2 of \cite{Burago}.\hfill $\Box$
%As an immediate consequence we get theorem \ref{min_thm}. 

\noindent
\section{The Hessian} {\sc Proof of Proposition \ref{Hessian}.}
\label{s:Hessian}

We continue with the same notation used in the proof of Proposition \ref{critical_pts} as above (subsection \ref{ss:critical_pts}), except now we add the shorthand:
\begin{equation}
n_{i,j} = n(q_i, q_j) \quad ; \quad r_{i,j} = |q_i - q_j|.
\label{def:nij}
\end{equation}
Write $d^2 f$ for the Hessian of a function.  Returning to the 
function  $|x|$ on $E$, a   routine computation shows that 
 $$d^2 r = \frac{1}{r} | d x - \langle n , dx \rangle n |^2; \qquad n = x/r .$$
 An application of the chain rule now  shows that  
 $$d^2 r_{12} =   \frac{1}{r_{12}} |(dq_1 -dq_2) - \langle n_{12} , (dq_1-dq_2)  \rangle n_{12} \Big|^2  $$
which simply means that
$$(d^2 r_{12} )_{(q_1,q_2)} (\xi_1, \xi_2) = \frac{1}{r_{12}} 
\big|(\xi_1 - \xi_2) - \langle n_{12} , (\xi_1-  \xi_2)  \rangle n_{12} \big|^2  $$
as a quadratic form on $E \times E$.

It follows that 
$$d^2 S_{A, B} = \sum_{i=0} ^k  \frac{1}{r_{i,i+1} }
\big| (dq_i - dq_{i+1}) - \langle (dq_i - dq_{i+1}), n_{i,i+1} \rangle n_{i,i+1} \big| ^2,$$
provided we set $q_0 = A$, $q_{k+1} = B$, $dq_0 = 0, dq_{k+1} = 0$.  
In other words
\begin{eqnarray}
d^2 S_{A, B} & = &  \sum_{i=0} ^k  \frac{1}{r_{i,i+1} } \big| (\xi_i -  \langle \xi_i ,  n_{i, i+1} \rangle n_{i, i+1}) -  (\xi_{i+1} -  \langle \xi_{i+1} ,  n_{i, i+1}  \rangle n_{i, i+1}) \big|^2 
\end{eqnarray}
is the Hessian  $d^2 S_{A,B} (\xi, \xi)$, a quadratic form in  
  $\xi = (\xi_0, \xi_1, \ldots, \xi_k, \xi_{k+1})$ where 
  we set $\xi_0 = \xi_{k+1} = 0$,  and where $\xi_i \in L_i$, $1 \le i \le k$.

Expand out this Hessian, focussing on the block diagonal terms:
$$d^2 S_{A,B} = \sum \frac{| \xi_i -  \langle \xi_i ,  n_{i, i+1} \rangle n_{i, i+1} |^2}{r_{i,i+1}}  + \frac{| \xi_i -  \langle \xi_i ,  n_{i-1, i} \rangle n_{i-1, i} |^2}{r_{i-1,i}}    + \text{ off-diagonal blocks}.$$
Since $A q_1 \ldots q_k B$ is 
  a point billiard trajectory the  projections of $n_{i-1, i}$ and $n_{i, i+1}$ onto $L_i$ are equal and so  
$$ | \xi_i -  \langle \xi_i ,  n_{i, i+1} \rangle n_{i, i+1} |^2 =  | \xi_i -  \langle \xi_i ,  n_{i-1, i} \rangle n_{i-1, i} |^2 = |\xi_i|^2 - \langle \xi_i, a_i \rangle ^2$$
where  we have written:   
$$\pi_i  (n_{i-1,i} ) = \pi_i  (n_{i, i+1})  : = a_i \, .$$
{\it We define this common  quadratic form on } $L_i$ to be:  
$$ \| \xi_i \|_i ^2 = | \xi_i -  \langle \xi_i ,  n_{i, i+1} \rangle n_{i, i+1} |^2 = |\xi_i|^2 - \langle \xi_i, a_i \rangle ^2$$
{\it and we observe that it defines a new inner product on} $L_i$.
Indeed since   $|a_i| < 1$ (equivalently  $n_{i, i+1}$, $n_{i-1, i}$ are unit vectors and are not in $L_i$), this quadratic form is indeed that  of a positive definite inner product on $L_i$.
Setting:  
$$\beta_i = \frac{1}{r_{i-1,i}} + \frac{1}{r_{i,i+1}} $$
we see that
$$d^2 S_{A,B} = \Sigma \beta_i \| \xi_i \|_i ^2 \ +   \text{ off-diagonal}.$$

We proceed to understand the off-diagonal terms. 
% In order to do so, set
After polarizing the quadratic form $d^2 S_{A,B}$ to obtain the associated symmetric bilinear form, still denoted $d^2 S_{A,B}$
we find that the off diagonal blocks are expressed in terms of the bilinear forms
$$Q_{ij} (\xi_i, \zeta_j) = \big\langle \xi_i -  \langle \xi_i ,  n_{i, j} \rangle n_{i ,j}\ ,\ \zeta_j -  \langle \zeta_j , n_{i,j}  \rangle n_{i,j} \big\rangle, \quad\text{with}\quad |i-j|=1 $$
with $\xi_i \in L_i$, $\zeta_j \in L_j$ and
%$$n_{ij} = \frac{p_i - p_j}{r_{ij}},  r_{ij} = |p_i -p_j |$$
so $Q_{ij}$ is an ``off-diagonal'' bilinear form:  
$$Q_{ij}: L_i \times L_j \to \R .$$
Then the off-diagonal terms of the polarized Hessian are: 
$$\text{off-diagonal terms}  =  -\sum_{|i-j| =1}  \frac{1}{r_{i  j}} Q_{ij} (\xi_i, \zeta_j) . $$
Now, using our $\langle \cdot, \cdot \rangle_i$ inner products we have that
$$|Q_{ij} (\xi_i, \zeta_j) | \le \| \xi_i \|_i \, \| \zeta_j \|_j$$
according to the usual Cauchy-Schwartz inequality on $E$.
It follows that if we define the operators $S_{ij}: L_j \to L_i$
by 
$$Q_{ij} (\xi_i , \zeta_j) = \langle \xi_i, S_{ij} \zeta_j \rangle_i\,,$$
then the  operator norms of the $S_{ij}$ are
\begin{equation} 
\| S_{ij}  \| \le 1
\label{opnorm}
\end{equation}
relative to the norms $\| \cdot \|_i$,  $\| \cdot \|_j$.

Endow  $\Lambda = L_1 \times L_2 \times \ldots L_k$ with  the inner product  $\langle \cdot , \cdot \rangle_*$
whose   squared norm is $\Sigma \|\xi_i \|_i ^2$.
Then we can define a $\langle \cdot , \cdot \rangle_*$--symmetric matrix $M: \Lambda \to \Lambda$
in the usual way:  
$$d^2 S (\xi, \zeta) = \langle \xi,  M \zeta \rangle_* $$
and we find that $M$ is block-tridiagonal with form:

\begin{equation} M  = 
\left(
\begin{array}{cccccc}
 \beta_1   &   -  \frac{1}{r_{12}}  S_{12} & 0 & 0 & \cdots & 0  \\
  -  \frac{1}{r_{21}}  S_{21 }&  \beta_2  &  -  \frac{1}{r_{23}}  S_{23}  & 0 &  \cdots & 0   \\
   0  &  -  \frac{1}{r_{23}}  S_{32}    & \beta_3 &  - \frac{1}{r_{34}}S_{34}     & 0   &  \cdots    \\
    \vdots &  \vdots &  \vdots &  \vdots &  \vdots  \\
      0  &   0 &  \cdots &  0 & - \frac{1}{r_{k-1,k}}  S_{k,k-1}   & \beta_k
\end{array}
\right)  \ .
\label{M}
\end{equation} 

In what follows, it is crucial to observe  that the coefficients  in all the rows but the first and last
satisfy a simple linear condition:   
$$ \frac{1}{r_{i,i-1}} +   \frac{1}{r_{i, i+1}} =    \beta_i .$$ 

In order to establish that $M$ is invertible and hence $d^2 S$ is nondegenerate
we form
$$P = DM$$
where $D$ is the block-diagonal matrix whose $i$th block is  $\frac{1}{\beta_i}$. Thus $D$
is the matrix for the   invertible transformation $(D \xi)_i = \frac{1}{\beta_i} \xi_i$. 
Proposition \ref{Hessian}  will be established once we establish 
the following lemma \ref{lem:invertible}. 
\hfill $\Box$
\begin{lemma} \label{lem:invertible}
$P$ is invertible.
\end{lemma}
\noindent
{\sc Proof.} 
We compute that
$$P = I -A$$
where $I$ is the identity and $A$ is tridiagonal block matrix with $0$'s on the diagonal and
\begin{equation} A  = 
\left(
\begin{array}{cccccc}
0   &   b_1  S_{12} & 0 & 0 & \cdots & 0  \\
 a_2  S_{21 }& 0  &  b_2 S_{23}  & 0 &  \cdots & 0   \\
   0  & a_3 S_{32}    &0  &  b_3 S_{34}     & 0   &  \cdots    \\
    \vdots &  \vdots &  \vdots &  \vdots &  \vdots  \\
      0  &   0 &  \cdots &  0 & a_k S_{k,k-1}   & {0} 
\end{array}
\right)
\label{A}
\end{equation} 
where 
$a_i = \frac{r_{i,i+1}}{r_{i,i-1}+ r_{i,i+1}}$, $b_i = \frac{r_{i,i-1}}{r_{i,i-1}+ r_{i,i+1}}$
so that
$$0 < a_i, b_i < 1, \text{ and }  a_i + b_i = 1.$$
To prove  that $P$ is invertible is equivalent to proving  lemma \ref{lem:not:eigen} immediately below.
\hfill $\Box$
\begin{lemma}  \label{lem:not:eigen}
$1$ is not an eigenvalue of $A$.
\end{lemma}
{\sc remark} The argument underlying lemma \ref{lem:not:eigen} 
was inspired by playing with  the situation in which all the $S_{i,j} = 1$
so that $A$ becomes  a tri-diagonal matrix with $0$'s on the diagonal, all  other tridiagonal entries   positive,
and all of its row sums except the first and last being $1$, that is,  a Perron-Frobenius  matrix.

\vskip .2 cm  
\noindent
{\sc Proof of lemma \ref{lem:not:eigen}.} 
Introduce the new norm on $\oplus L_i$:
$$\| \xi \|^* : = \max_j \| \xi_j \|_j\, .$$
Suppose that
$A \xi = \xi$.
We must show that $\xi = 0$.
The eigenvalue equation reads:
$$\xi_i = a_i S_{i, i-1} \xi_{i-1} + b_i S_{i, i+1} \xi_{i+1} ,\quad 1 < i < k$$
together with
$$\xi_1 = b_1 S_{1,2} \xi_2
\quad\mbox{and}\quad
\xi_k = a_{k} S_{k,k-1} \xi_{k-1}.$$
By way of contradiction, suppose that $\xi \ne 0$
so that $\| \xi \|^* > 0$. 
Let $i$ be an index such that 
$$\|\xi\|_i = \max_j \| \xi_j \|_j := \|\xi \|^* \, .$$
Then $i$ cannot be $1$ or $k$, for if it were,
taking norms of the eigenvalue equation for these indices and use eq (\ref{opnorm}),
together with  $|b_1|, |a_k| < 1$
we would have
$\| \xi \|^* < \|\xi \|^*$, a contradiction.
So $1 < i < k$.
Taking   norms we get
$$\|\xi_i \|_i    \le a_i  \|S_{i, i-1} \xi_{i-1}\|_i + b_i \|S_{i, i+1} \xi_{i+1}\|_i  
\le a_i \|\xi_{i-1}\|_{i-1} +  b_i \|\xi_{i-1}\|_{i-1}\, .$$
Now $\|\xi_{i\pm1}\|_{i \pm 1}  \le \| \xi \|_i$ and $a_i + b_i = 1$.  It follows that unless both
$\|\xi_{i-1}\|_{i-1}$ and $\|\xi_{i+1}\|_{i+1}$ are equal to $\|\xi_i \|_i$ we will have again that
$\| \xi \|^* < \| \xi \|^*$, a contradiction.
Thus we have now have that 
$$\| \xi \|^* = \| \xi_i \|_i = \|\xi_i \|_{i-1} = \|\xi_i \|_{i+1}.$$
Continuing in this manner we march up or down the indices until eventually
$\|\xi_1 \|_1 = \|\xi_k \|_k = \|\xi \|^*$
and we return to our original contradiction. 

Finally that the Hessian is positive definite and not simply nondegenerate follows
from theorem \ref{min_thm}, and hence the CAT(0) ideas of \ref{sect:CAT}.  \hfill $\Box$

 \section{Thickening}
 \label{sec:thick}
 
 In this section we  thicken each  subspace $L \in \mL$ and in so doing
 obtain an approximating  deterministic dynamics to our point billiard system.     We  introduce the notion of  trajectories
 being transverse.  We prove that  every transverse  point billiard solution in $\mB (L_1 \ldots L_k)$   is the limit of a family
 of thickened billiard trajectories as the   thickening parameter
 tends to zero.  
   This limit assertion yields another perspective on point
  billiards as well as a new proof of the main parts of theorems~\ref{main} and \ref{main_unreduced}.

 \begin{definition}  \label{def:r:thickened} Choose  positive
scale factors $\sigma_L >0$ for each $L \in \mL$.  
For $r>0$ ,   $L \in \mL$ set  
\[L^{(r)}:=\{q\in E\mid  d(q,L)\le \sigma_L r\}, \]
\[Z^{(r)}  :=\{q\in E\mid  d(q,L)=  \sigma_L r\} ,  \qquad\mbox{and} \]
\[\qquad {\mathcal M}^{(r)}:={\rm cl}\big( E - \textstyle \bigcup_{L\in {\mathcal L}} L^{(r)} \big) .\]
An ``r-thickened billiard trajectory'' is a solution $q^{(r)}$ to the deterministic billiard problem played on the table  ${\mathcal M}^{(r)}$.
\end{definition}

 The walls of our table ${\mathcal M}^{(r)}$ are the unions of the cylindrical hypersurfaces   $Z^{(r)}$ minus certain small `corner' or   intersection   parts
 where two or more of the interiors $L^{(r)}$ of these cylinders   intersect. Away from these small corners,
 the $r$-thickened billiard problem is  a deterministic dynamics of standard billiard type.

\begin{example}[Thickened $N$-body billiards = ideal gas] The  reason behind   the   scale factors in  definition \ref{def:r:thickened} arises here.    
 The formula
 \begin{equation} \sigma_{ab}  r_{ab} = {\rm dist}_E (q, \Delta_{ab}); \qquad \sigma_{ab} =  \sqrt{\frac{m_a + m_b}{m_a m_b }}.
 \end{equation}
 relates the  usual  
distance $r_{ab} = |q_a -q_b|$   between the $a$th and $b$th  bodies and the  $E$-distance   
between the corresponding  configuration point $q = (q_1, \ldots , q_N) \in E = (\R^d)^N$ and the collision subspace $\Delta_{ab}$
 (See,  for example, the proof of lemma 2 and eq (4.3.15a) in \cite{Mont}.)  
If we  take  $ \sigma_{ab}$ for each $L = \Delta_{ab}$  in definition \ref{def:r:thickened} 
 then the  domain ${\mathcal M}^{(r)}$ within which the thickened billiard moves is precisely the configuration space  of $N$ hard balls with centers $q_a$ and radii  $r/2$,
 that is, an ideal gas (but unconfined to a box).  
 %The deterministic billiard dynamics is that of the $N$ hard balls of  diameters all $r$ and masses $m_a$. 
\end{example} 

\begin{definition} 
If $q$ is a point billiard solution then an  {\bf $r$--family for} $q$ 
is a family of  r-thickened billiard trajectories  $q^{(r)}:\R\to {\mathcal M}^{(r)}$,
$r \le r_0$, some $r_0 > 0$, such that 
$q^{(r)} \to q$ in the compact-open topology as $r \to 0$.
\end{definition}

We would like to say that every point billiard trajectory admits an $r$-family.
But that is not true.  However the exceptional trajectories  are quite easy to understand.
Our reflection rule (eq \ref{eq:momentum}) allows for  point billiard trajectories which pass right through a collision subspace
without changing direction: $n_{i-1,i} =n_{i,i+1}$.  For deterministic  billiards  this  cannot happen: collisions with walls
change  direction.
% (ignoring tangencies!). 

\begin{definition}$\quad$\\ \label{def:transverse}
An internal vertex of a polygonal path  $q_0 q_1  \ldots q_k q_{k+1} $ is a vertex $q_i$
such that the  edges $q_{i-1} q_i$ and $q_i q_{i+1}$ incident to it  form a line segment $q_{i-1} q_{i+1}$
with $q_i$ in the interior.  A polygonal path is {\bf transverse} if it has no internal vertices.
\end{definition}

%\begin{definition}$\quad$\\ \label{def:transverse}
%A point $(q_0,q_1, \ldots,q_k,q_{k+1})\in E^0\times \Lambda^0\times E^0$
%is called {\bf transverse}  if the velocities (see \eqref{def:nij})
%before and  after  each  collision are different: $n_{i-1,i}\neq n_{i,i+1}$ for $i=1,\ldots,k$.\\
%The {\bf transverse set} is the subset of $E \times \Lambda \times E$ 
%consisting of all transverse points.  
%\end{definition}

Here is the  main result of this section. 
\begin{proposition}  
\label{lem:r-thickened-shadows}
Any transverse point billiard trajectory $q$ admits  an $r$--family $q^{(r)}$. 
\end{proposition}
{\sc Basic remark.} The  itineraries of each path in the $r$--family agree with those of their limit  $q$, 
 for all   $r$  small enough.  

{\sc Caveat.} The set of transverse billiard trajectories need not be   dense within the set of all
billiard trajectories realizing a given itinerary. 
\begin{example}
If   $L\in {\mathcal L}$ has codimension one  then `half' of  the  transverse billiard trajectories 
colliding with $L$  are scattered back into the same half space of $E\setminus L$
while the other half  pass straight  thru $L$ without their direction
of travel being altered. So the transverse $L$-colliding trajectories are not transverse. 
\end{example}
\begin{example}
If three consecutive different  scattering subspaces $L_{i-1}$, $L_i$ and $L_{i+1}$ are coplanar lines
then any trajectory which has  $L_{i-1} L_i L_{i+1}$ as part of its itinerary will  have
$n_{i-1,i} = n_{i,i+1}$ and hence is not transverse. 
\end{example}

Our proof of  proposition \ref{lem:r-thickened-shadows}   relies on a minimizing property of thickened billiard trajectories
quite similar to that  used for our  earlier point billiards arguments 
but with one crucial difference.  The difference is the existence of  ``ghost billiards''. (See lemma \ref{lem:ghost}.)

Thicken our old parameter space  to 
%  the thickened   parameter space
$$\Lambda ^{(r)} = L_1 ^{(r)} \times \cdots \times L_i ^{(r)} \times \cdots \times L_k^{(r)}. $$
and consider the polygonal path length function  with $\Lambda ^{(r)}$
as the input vertices to form the thickened analogue of 
$S_{AB}$.

%We call such a trajectory ``transverse'' if $n_{i-1,i} ^r\neq n_{i,i+1} ^r$ where we hope the
%meaning of the superscript notation is obvious here.  
\begin{lemma} 
\label{lem:r-min}
 For fixed $A, B \in int({\mathcal M}^{(r)}) $ the minimum of $S_{A,B}$
over $\lambda^{(r)} \in \Lambda ^{(r)}$ exists and is unique. Moreover, any
local minimizer or critical point for $S_{AB}$  is this global minimizer. 
 If that global   minimizer is   transverse then it  is a solution to the  deterministic billiard problem in $int({\mathcal M}^{(r)})$.
 Conversely, any solution to the deterministic billiard problem is a minimum
of $S_{AB}$ where $A, B$ are taken on the incoming and outgoing rays of the solution.
\end{lemma}
\begin{lemma}[Ghost Billiards]
\label{lem:ghost} 
There exist non-transverse minimizers.   For these,  the interior vertex  $q_i ^{(r)}$ is part of a line segment  $q_{i-1}  q_{i+1} ^{(r)}$
which is   either tangent to  $Z_i ^r$ at $q_i ^{(r)}$  or (the more important case) 
which passes through the   interior of  $L_i ^{(r)}$ so that $q_i ^{(r)}$ may be taken to lie in that interior,  in which case we  say that the minimizer is a 
    ``ghost billiard trajectory'' in honor of what such a trajectory  looks like in the   thickened $N$-body billiards case. 
\end{lemma}  

{\sc Proof of  lemma \ref{lem:r-min}.}  With the exception of the assertion regarding transverse minimizers,  the proof of lemma \ref{lem:r-min} is almost identical to the
proof of the minimization property which we gave above in propositions \ref{critical_pts} and \ref{Hessian}    for    point billiard trajectories.  A
proof can  also be found in \cite{BFK1}.   
The unique  global minimization property  of $S_{AB}$ is  achieved in  a manner identical to our proof
of theorem \ref{min_thm}.  We form  Hadamard spaces
by gluing sheets   $E_i = E$  together, now along the convex bodies $L_i ^{(r)}$.
To understand the assertion regarding transverse minimizers within the lemma we must  understand a bit about
the ghost billiards of  lemma  \ref{lem:ghost}.  

{\sc Sketch, Proof of  lemma \ref{lem:ghost}.} Take the case of a single $L$, that is, of an itinerary of length 1.  Set $K =  L^{(r)}$,   a convex body with non-empty interior in $E$.
 Suppose that the line segment $AB$ passes through the interior of $K$.   Put $A \in E_0$
 and $B \in E_1$ in the gluing construction $E_0  \cup_K E_1$ of Reshetnyak's theorem.  The geodesic from $A$ to $B$
 is now the  straight line segment $AB$  passing through the convex body without being deflected and still passing from one  sheet to the other.
{\it This is our ghost geodesic!  } Take $q_1 \in AB \cap K$  when  minimizing the thickened  $S_{AB}$ to arrive at the  non-transverse  minimizer $A q_1 B$.  If, on the other hand, a minimizer is transverse
it cannot be a ghost billiard (nor can it be a billiard with a tangency to $K$).  Such a transverse minimizer must    correspond to an ``honest billiard'' - a solution
to the deterministic billiard system.

{\sc Proof of Proposition \ref{lem:r-thickened-shadows}}.
Let $q$ be a transverse point   billiard trajectory with vertices $q_i \in L_i$ listed in order. 
Choose points $q_0 = A$ and $q_{k+1} = B$ on the ingoing and outgoing rays. 
% as per XXXW and perform a time translation so that $q(0) = A$.  
By a slight abuse of notation, we will also write $q$ for that finite  part of $q$
joining $A$ to $B$.    Let   ${\mathcal P}_0$ 
denote  the space of polygonal paths $q'$ starting at $A$,
ending at $B$ and having $k$ vertices $q'_i \in L_i$ in between, listed in order.   Then $q \in {\mathcal P}_0$.  According to the set-up 
from the  beginning of section \ref{s:generating_families}, 
${\mathcal P}_0$ is naturally isomorphic to the parameter space $\Lambda$ and 
 $S_{AB}$ is the restriction of the length functional $\ell$
to ${\mathcal P}_0$.     Proposition \ref{critical_pts} and theorem  \ref{min_thm}  assert that    $q$ is the  global minimum of $S_{AB}$.
Write $T = S_{AB} (q)$ for this minimum value.  (Thus, since $q$ is parameterized by arclength
if $q(0) = A$ we have that  $q(T) = B$.)

{\it Claim 1.}  There is a small  positive constant  $\delta_0$  and a positive   constant $K$ 
with the following significance.  If $\delta < \delta_0$ 
and if    $q' \in {\mathcal P}_0$ has the property that one of its  $k$ corners $q' _i$ 
satisfies   $|q'_i - q_i| \ge  \delta$ then  $\ell(q') \ge  T + K \delta ^2$.

\vskip .2cm

Write ${\mathcal P}(r)$ for the space of polygonal paths starting at $A$, ending at $B$
and having $k$ vertices $q''_i \in Z_i (r) = \partial L_i ^{(r)}$. 
For  $q'' \in {\mathcal P}(r)$ we write $\pi(q'') \in {\mathcal P}_0$
for the polygonal path whose $k$ vertices are $q'_i = \pi_i (q''_i) \in L_i$ where  $\pi_i: E \to L_i$ is the orthogonal projection.

{\it Claim 2.}
Suppose that   $\delta, \delta_0$  and $K$  are as in Claim 1.  Take    $r_0 >0 $ such that   $2 \Sigma \sigma_i r_0 < K \delta^2 /2$
where $\sigma_i = \sigma_{L_i}$ are the scale factors  attached to $L_i$ as per definition~\ref{def:r:thickened}.  If   $r < r_0$ and 
 $q'' \in {\mathcal P}(r)$ is such that  $\pi(q'') \in  {\mathcal P}_0 $ satisfies  the hypothesis of claim~1 (i.e.   $|q' _i -q_i  | \ge  \delta$ for some $i$)
 then  $\ell(q'') > T + K \delta ^2 / 2$.

{\it Claim 3.} [Curve Shortening]  For any $q \in {\mathcal P}_0$ and any $r$ sufficiently small  there is
a  $q'' \in {\mathcal P}(r)$ with $\ell(q'') < \ell(q)$.

We now show how the three claims yield the lemma.  Afterwards we prove the claims. 
 By lemma \ref{lem:r-min}   
for all $r$ there exists a unique length  minimizer $q_* = q_* ^{(r)} \in {\mathcal P}(r)$.
Apply the curve shortening process of Claim 3 to $q$ in order to  get a $q'' \in {\mathcal P}(r)$
with $\ell(q'')   < T = \ell(q)$.  Thus $\ell(q_*) \le \ell(q'') < \ell (q) =  T < T + (K/2) \delta^2$.
Take $\delta, \delta_0, r_0, r$ as per claim 2 
 to conclude that  each  vertex $q' _i$ of the  projected polygonal curve $q' = \pi(q_*) \in {\mathcal P}_0$
   satisfies   $|q'_i - q_i| < \delta$.  But $|q_{* i} - q_i| = \sqrt{|q_{*i} - q'_i|^2 + |q'_i - q_i|^2 } = \sqrt{ \sigma_i ^2 r^2 + |q'_i - q_i|^2} < \sqrt{2} \delta$
   where in the first equality we used the fact that $q'_i = \pi_i (q_{*i})$ is the orthogonal projection onto $L_i$. 
    Letting $\delta \to 0$ we get that $q_{*i} \to q_i$.   Since $q$ is transverse, eventually, for $r$ small enough,  $q_*$ is also transverse,
    and hence a thickened billiard solution.  This proves that the $q_*$ form an $r$-family for $q$.    
   
   It remains to prove the three claims. 
    
    {\it Proof of  Claim 1.}  (i) Since the Hessian of $S_{AB}$ is positive definite at $q$ (proposition \ref{Hessian})
    we have that there exists a $\delta_1 >0 $ such that whenever $q' \in {\mathcal P}_0$ satisfies
    $|q' _i - q_i| < \delta_1$  for all $i$ then $S_{AB} (q') - S_{AB} (q) \ge \Sigma_i  K |q' _i - q_i|^2$.
    Our $\delta_0$ will eventually be less than or equal to $\delta_1$. 
    
    (ii) Now if $q'  \in {\mathcal P}_0$ has any vertex $q'_i$ such that    $|q'_i -q_0| \ge T + K \delta_1 ^2$ then $\ell(q') >   T + K \delta_1 ^2$.
    
   (iii) Now  restrict the length functional  to the compact set of paths $q' \in {\mathcal P}_0$  all of whose vertices $q'_i$ satisfy  $|q' _i - q_i| \le  T+ K \delta_1 ^2$
   and at least one of which satisfies $|q' _i - q_i | \ge  \delta_1$.  This set of  polygonal    paths 
   is naturally  homeomorphic to a  compact subset of $\Lambda$, and  as such the   length functional  $S_{AB}$ achieves its minimum value $T_M = \ell(q_M)$ on the
    set.    $T_M  > \ell(q) = T$   because $q$ is the unique global    minimizer of $S_{AB}$ and $q$ is not in our compact set.
    Write $T_M = T + \epsilon_M$ to get that   $\ell(q') > T + \epsilon_M$ for all the paths in this compact set. (We have   $\epsilon_M \le \delta_1$.)
  
  Combining (iii) and (ii) we  see that  if any 
    path $q' \in {\mathcal P}_0$ has one vertex $q'_i$ with $|q' _i - q_i | \ge  \delta_1$  then  $\ell(q') \ge T + \min\{ \epsilon_M,  K \delta_1 ^2 \}$.
    Choose $\delta_0$ so that $\min\{ \epsilon_M,  K \delta_1 ^2 \} = K \delta_0 ^2$.   This $\delta_0$
    will do the needed  trick. For let $\delta \le \delta_0$ and suppose that $q'$ has one vertex $q'_i$ with $|q'_i - q_i| \ge \delta$. 
    Let  $i$ be an index such that $|q'_i - q_i|$ is maximized  and let this maximum value be $m$.  Thus $m \ge \delta$.   If $m \ge \delta_1$, then
        from the previous paragraph   $\ell(q') \ge  T +  \min\{ \epsilon_M, K \delta^2 \} =  T + K \delta ^2$.
    Otherwise, $m < \delta_1$ and the   Hessian bound holds on $q'$, yielding  
    $\ell(q') \ge T+  K \Sigma_i |q'_i - q_i|^2 \ge T + K m^2 \ge   T+ K \delta^2$. 
    $\Box$.

        {\it Proof of  Claim 2.}
        
    Suppose that $q''$ is as in the statement of this claim so that  its projection $q'$ satisfies the conditions
    of Claim 1 and thus  $\ell(q') \ge T + K \delta^2$.
    The difference  between the  vertices of $q''$ and $q'$ satisfies 
    $|q''_i - q'_i| = \sigma_i r^2$ since the projection $\pi_i$ is orthogonal
    and $q''_i \in Z_i (r)$.   By the  triangle inequality
$|q'_i - q'_{i+1}| - \sigma_i r - \sigma_{i+1} r \le |q''_i - q''_{i+1}|$
so that
$$T + K \delta^2 \le \ell(q') - 2 \Sigma \sigma_i r \le \ell(q'').$$
By assumption $2 \Sigma \sigma_i r \le (K/2) \delta^2$, yielding the desired result, 
$$ T + (K/2) \delta ^2   \le \ell(q'')$$

\begin{figure}
\begin{center}
\includegraphics[height=10cm]{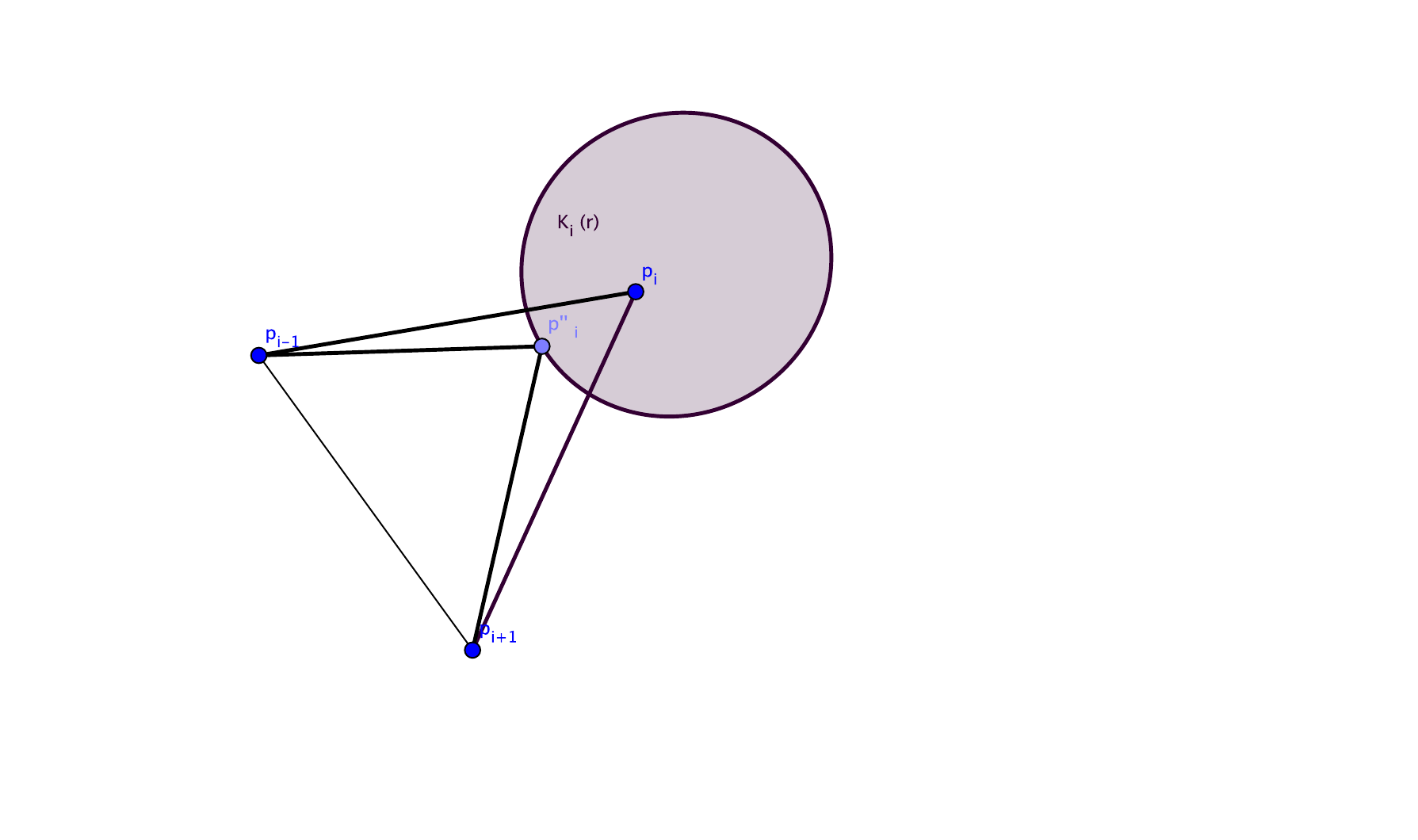}
\vspace{-5mm}
\end{center}
\caption{The curve shortening process of Claim 3}
\label{fig:curveshortening}
\end{figure}

{\sc Proof of Claim 3.} For each vertex $q_i$ consider the  triangle $\Delta_i$ whose vertices  are
$q_{i-1}, q_i, q_{i+1}$.  By the transversality condition $\Delta_i$ is a nondegenerate
triangle and so   lies in a unique affine planes  $P_i \subset E$.  
The  solid cylinder  $L_i (r)$ intersects $P_i$ in a convex domain $K_i (r)$  (the interior of an ellipse) containing  the vertex $q_i$, 
and for $r$  sufficiently small  the other two vertices of our triangle are not in $K_i (r)$.   See figure \ref{fig:curveshortening}.  
$K_i (r)  \cap \Delta_i \subset P_i$  is a convex planar domain whose   boundary  consists of three arcs, two being line segments forming part of the  edges of $\Delta_i$
and the   third  curved arc  $C_i$  (a subarc of the ellipse) lying in the interior of $\Delta_i$. 
Choose any point $q''_i  \in C_i$ on this third arc. Then $q''_i \in Z_i (r)$.  But  $q''_i$ is in the interior of our triangle, and so by a general property of interior points of triangles we have 
that  $|q_{i-1} - q''_i| + |q''_i - q_{i+1}| < |q_{i-1} - q_i| + |q_i - q_{i+1}|$.  See figure \ref{fig:curveshortening} again. 

It follows that the replacement of vertex $q_i$ by $q''_i$  leaving all other vertices of $q$ unchanged 
shortens the polygonal path.  Indeed,  only the edge lengths of the two edges 
incident to the changed vertex change and the sum of these two decrease. 
Write this replacement process as $q \mapsto \sigma_i (q)$. 
Apply  this polygonal curve shortening   process consecutively,  
vertex-by-vertex, 
$q \mapsto \sigma_1 (q) \mapsto \sigma _2 (\sigma_1 (q)) \mapsto \ldots \mapsto \sigma_k ( \ldots (\sigma_1 (q)) \ldots )) =q''$.
Each replacement shortens the path. The  $k$th application yields a 
path $q'' \in {\mathcal P} (r)$ all of whose corners $q'' _i$ lie on their respective boundary sets   $Z_i (r)$ and which is shorter than the original path $q$.
$\Box$ (for the proof of Proposition~\ref{lem:r-thickened-shadows}). 

\subsection{Relation to proofs of theorems \ref{main}, \ref{main_unreduced}}

The $r$-thickened dynamics is deterministic and symplectic. The graph of its time $T$ flow is,
morally speaking, a   Lagrangian graph.   This graph is partitioned  up into pieces according to the itineraries of the
trajectories.  In some sense, which we are purposely  vague on, our Lagrangian relations of theorems
\ref{main} and \ref{main_unreduced}, are the limits $r \to 0$ and $T \to \infty$ of the pieces of this graph.
Among the problems faced in turning this idea into a complete proof is the fact that the flow is
not continuous due to the instantaneous velocity changes suffered at collisions.

 \section{Conservation  Laws, Symmetries, and Scaling.}  
 \label{s:symmetries} 
 
 Solutions to the $N$-body problem enjoy conservation of linear and angular momentum.
We expect that our $N$-body  billiard trajectories to 
 obey  these same conservation laws.   They do.  We show
derive the laws  from  the   group invariance   of the   collision subspaces.  
 We end this section with a remark on a scaling symmetry for billiard trajectories  and what it implis for the closures of our Lagrangian relations
 in theorem \ref{main}. 
 
 \subsection{Momentum maps for free motion and its restrictions.}
 
  The usual  linear and angular momentum are  the  components of the momentum map for the
  group  of rigid motions of the underlying Euclidean space.  We take, to start with,
  the full group of rigid motions of our $E$, and later restrict to subgroups mapping the collision
  spaces to themselves.    
  
  The  group  ${\rm Iso}(E)$ of rigid motions of $E$ splits into translations and rotations. 
  Write ${\rm Lie}({\rm Iso}(E))^* $ for the dual of the Lie algebra of this Lie group.  Using
  the inner product on $E$ we have canonical identifications:  $T^* E = E \times E$
 and ${\rm Lie}({\rm Iso}(E))^* = {\rm Lie}({\rm Iso}(E)) = E \oplus \Lambda ^2 E$.  The full  momentum map is  then  
 the  map $\Phi = (P, J) : E \times E \to E \oplus \Lambda ^2 E$ whose  
 first (translational)  factor $P(x,v) = v$ we   call  the ``full''  linear momentum
 and whose   second (rotational) factor    $J_{\rm rot} (x, v) = x \wedge v$ 
 we call  the `full' angular momentum and is that for the full rotation group.
 
 Free (=straight line) motion $(x,v) \mapsto  (x + tv, v)$ is a Hamiltonian flow on
 $E \times E$ which has the full momenta as conserved quantities.
  
Now restrict attention to the  Lie subgroup    of   those $g \in {\rm Iso}(E)$
 such that $g(L) = L$ for each $L \in \mL$.  
Being a subgroup of ${\rm Iso}_+ (E)$,  this subgroup also  acts symplectically  on   the phase space $E \times E$
 and  has its own  momentum map which is well-known to be the composition of the previous full momentum map  $\Phi$
  with the orthogonal projection onto our subgroup's Lie algebra.
   In this way we get   linear and angular momenta associated to our collision-preserving  subgroup. :
 $$P_{\rm tr} (x,v) = \pi_{\rm tr} (v) \in L_{\rm tr}$$
 and
 $$J (x, v) = \pi (x \wedge v) \in {\rm Lie}(H)$$
 where $\pi_{\rm tr}: E \to L_{\rm tr}  $ projects onto  the translational part of our subgroup
 and $\pi:  \Lambda ^2 E \to {\rm Lie}(H)$
 projects onto its rotational part.  In the next two subsections we compute these
 projections and derive their conservation consequences.

 %The rotational part $\{g \in O(E):  g(L) = L \text{ for all } L \in \mL \}$ 
% will be generally harder to compute explicitly.     
 
 \subsection{Linear Momentum and Translation invariance}  
 
  The translational part of our collision-preserving subgroup is:  
  \begin{equation}
 \label{eq:translation_subspace}
 L_{\rm tr} = \bigcap_{L \in {\mathcal L}} L.  
 \end{equation}
  In other words, $L_{\rm tr}$  is precisely the subgroup of translations of $E$ which maps {\it each}  $L$ onto itself.
  Write $\pi_{\rm tr}: E \to L_{\rm tr}$ for the orthogonal projection
  onto this subspace, as above, we have:  
  \begin{proposition} The `total linear momentum'' $\pi_{\rm tr} (v)$
  is constant along each billiard trajectory.
  \end{proposition}
  
  {\sc Proof. } At each collision we have $\pi_L (v_-) = \pi_L (v_+)$.
  But $L_{\rm tr} \subset L$ for all $L \in \mL$.
  So $\pi_{\rm tr} (v_-) = \pi_{\rm tr} (v_+)$ at each collision: 
  the total momentum remains unchanged at each collision
  and thus $\pi_{\rm tr} (v)$ is constant along any given billiard trajectory.
  \hfill$\Box$
  
  {\sc $N$-body billiard momentum conservation }  In $N$-body billiards the   intersection of all of the $\Delta_{ij}$ consists  
  of the $d$-dimensional subspace consisting of all vectors of the form $(z, z, \ldots, z)$, $z \in \R^d$.  
  It is the   subspace of $E = (\R^d) ^N$ generated by the translation group of $\R^d$. 
 The projection of a velocity $v \in (\R^d) ^N$ onto this subspace, relative to the mass metric, 
 is $(v_1, \ldots , v_N) \mapsto \Sigma m_a v_a$ which co-incides with total 
  linear momentum.  
 
 \subsection{ Angular Momentum and  Rotational Invariance}
 
 Now we consider the rotational part of our collision-preserving subgroup.  
 Denote this subgroup as  $H \subset {\rm O}(E)$ so that $H$ 
 consists of all rotations which map the collision subspaces to themselves.   
 We write $\pi_H : \Lambda ^2 E \to {\rm Lie}(H)$ for the orthogonal projection,
 identifying ${\rm Lie}(H)$ with a linear subspace of $\Lambda^2 (E)$.  Use the 
 naturally induced  invariant   inner product on $\Lambda^2 (E)$. On bivectors   $v \wedge w$
 the squared length for this inner product is  
 \begin{equation}
 \langle v \wedge w, v \wedge w \rangle = \det
\left(
\begin{array}{cc}
 \langle v , v \rangle  &  \langle w , v \rangle     \\
  \langle v , w \rangle &      \langle w , w \rangle
\end{array}
\right)
\label{eq:GramDet}
\end{equation}

  \begin{proposition} The `total angular momentum''  $\pi_{H} (x,v)$
  is constant along each billiard trajectory.
  \end{proposition}
  
  {\sc Proof.} Let $\xi \in {\rm Lie}(H)$.  Thus $e^{t \xi} \in H$ is a one-parameter family of rotations
  leaving each $L$ invariant.   The $\xi$-component of the full angular momentum is 
  $$J^{\xi} (x,v) = \langle x \wedge v,  \xi \rangle$$
  where $\langle \cdot , \cdot \rangle$
  is the inner product on   ${\rm so}(E) \cong \Lambda ^2 E$ just described. 
    Since $x \wedge v$
  is constant along straight line motions, $J^{\xi} (x,v)$ remains
  constant along the straight line segment parts of a billiard trajectory.    We must
  show that it remains unchanged at collisions. 
  The jump  in $J^{\xi}$ at a collision  with an $L$ at $q \in L$ is:  
  $$J^{\xi} (q,v_+) - J^{\xi} (q, v_-) = \langle q \wedge (v_+ - v_-), \xi \rangle$$ 
  Now from $\pi_L (v_+) = \pi_L (v_-)$ we have that $v_+ - v_- \in L^{\perp}$.
  Thus $q \wedge (v_+ - v_-) \in L \wedge L^{\perp} \subset \Lambda ^2 E$.  
  Our proposition now follows from the computational lemma:
  
  \begin{lemma} If  $\xi \in {\rm Lie}(H)$, $q \in L$ and $v \in L^{\perp}$
  then $\langle q \wedge v, \xi \rangle = 0$.
  \end{lemma}
  
 {\sc  Proof of lemma. }
  Using bilinearity of the inner product and formula (\ref{eq:GramDet}) one    verifies that for any $\xi \in \Lambda^2 E$
  we have  
  $$\langle q \wedge v, \xi \rangle = \langle \xi(v), q \rangle$$
  where on the right-hand side we view  are viewing $\xi$ as a skew symmetric map $\xi: E \to E$ using the canonical identification
  $\Lambda ^2 (E) \cong {\rm so}(E)$.   (Under this identification the bivector $q \wedge v$
  becomes the linear transformation  $e \mapsto (q \wedge v )(e) = \langle v, e \rangle q - \langle q, e \rangle v$ of $E$.)
  It follows that we also have
  $$\langle q \wedge v, \xi \rangle = -\langle \xi(q), v \rangle$$
  Now if $\xi \in {\rm Lie}(H)$ then   $e^{t \xi} \in H$ is a one-parameter family of rotations
  leaving each $L$ invariant.  Differentiating, we see that if   $q \in L$  then $\xi(q) \in L$. 
  But in the lemma $v \in  L^{\perp}$ so that  $ - \langle \xi (q), v \rangle = 0$.
  \hfill$\Box$
 
 {\sc $N$ body billiards.} 
 The group $H = {\rm O}(d)$ acts diagonally on the $N$-body configuration space $(\R^d)^N$
 leaving each $\Delta_{ij}$ invariant.   
The mass-metric projection of $q \wedge v \in {\rm so}((\R^d)^N) = {\rm so}(E)$ to 
 ${\rm Lie}(H) = \Lambda ^2 \R^d$ is   $\pi_H (q \wedge v) = \Sigma_{ a = 1} ^N m_a q_a \wedge v_a \in \Lambda^2 \R^d$,
 the usual formula for the total  angular momentum.  We have that total angular
 momentum is conserved for   $N$-body billiards.
 
 \subsection{Scaling and the Scattering map as a Legendrian Map.}   If $q(t)$ is a  billiard trajectory with itinerary $L_1 L_2 \ldots L_k$  then so is  
 $\lambda q (\frac{t}{\lambda}):= q_{\lambda} (t)$,  $\lambda > 0$.  
 Letting  $\lambda \to 0 $ brings all the collision points $q_i$ of $q_{\lambda} $  to the origin.  In this way we see  that 
  the closure of the  Lagrangian relation $\mathcal R$ for $\mB (L_1 L_2 \ldots L_k)$
 (see theorem  \ref{main}) 
 contains  points lying  in the Lagrangian relation for total collision described in  example \ref{ex:1} - namely  the product of the two zero sections of
 $T^* S(E) = \LINES(E)$.  
 
Scaling  acts on pairs $(q,v) \in E \times E_v$ by $\lambda (q,v) = (\lambda q, v)$. Directions are 
 left unchanged while 
 `impact parameters' $q$ are scaled.       
 %Fix an itinerary $L_1 L_2 \ldots L_k$  of length $k \ge 2$ and 
 Let $\tilde \mR = \tilde \mR (L_1 \ldots L_k)$
 be  the  ``unreduced'' Lagrangian relation of theorem \ref{main_unreduced}.  The scale invariance  of $\mB(L_1 L_2 \ldots L_k)$
 implies that $\tilde \mR$ is scale-invariant: $(q_A, v_A), (q_B, v_B) \in \tilde \mR \iff
 ((\lambda q_A, v_A),  (\lambda q_B, v_B) ) \in  \tilde \mR$. In other words, the
 Lagrangian relation is a scale-invariant submanifold of $(E \times E_v) \times (E \times E_v)$.
 
 Scaling commutes with reduction, and so induces a scaling action on $T^* S(E_v) \times T^* S(E_v)$
 which leaves the Lagrangian relation $\mR$ of theorem \ref{main} invariant.  In terms of
 the coordinates on $LINES(E) \cong T^* S(E_v)$  discussed in subsection \ref{sss:reduced},  scaling acts again by $(v, Q) \mapsto (v, \lambda Q)$.  
 Thus we expect to be able to form the quotient by this action to arrive at
a submanifold $\mR/ \R^+  \subset  (T^* S(E_v) \times T^* S(E_v) )/ \R^+$. 
 The latter is a nice manifold {\it provided we  delete its zero section} before forming the scaling quotient.    Indeed, for any manifold $X$,
 let $Z_X \subset T^* X$ denote the zero section of its cotangent bundle. 
 Then  $(T^*X \setminus Z_X)/\R^+ = \bP (T^* X)$   is a canonical contact manifold which fibers over
 $X$ with fibers $\R \bP ^{n-1}$'s,  $n = dim(X)$.  (See  \cite[Appendix 4]{Arnold}.) 
 We apply this observation  to $X = S(E_v) \times (S(E_v)$, using  
  $T^* (S(E_v) \times T^* S(E_v) = T^* (S(E_v) \times (S(E_v))$ to arrive at:  
 \begin{theorem} 
 \label{Legendrian}  If  the itinerary   has length  greater than 1 then 
 the quotient $\mR / \R^+$ of the      Lagrangian relation $\mR$ of  theorem \ref{main}  by the scaling group $\R^+$
 is a    submanifold of $\bP T^* (S(E_v) \times (S(E_v) )$ of dimension $2 dim (E) -3$
 which is Legendrian  relative to the 
 (nonstandard) contact form 
 $$\Theta = {\vec Q_-} \cdot d {\vec v_-} - {\vec Q+} \cdot d {\vec v_+}. \footnote{The form itself varies under scaling so is not well-defined as a one-form on
 $\bP T^* (S(E_v) \times (S(E_v) )$. The form is to be viewed projectively:  its zero locus, which is  the contact distribution, is independent of scaling.} $$ 
 %on $\bP T^* (S(E_v) \times (S(E_v) )$
 The  projections  to the incoming and outgoing velocity spheres are scale invariant maps and combine to  yield the  Legendrian fibration 
 $$\pi_- \times \pi_+ :\bP T^* (S(E_v) \times (S(E_v) ) \to S(E_v) \times S(E_v)$$
 under which  the image of  $\mR / \R^+$ is a  (possibly singular) hypersurface provided it is transverse to the fiber at some point
 % in the product of velocity spheres.  
 % is transverse to the fibers    then it projects  locally diffeomorphically onto an immersed hypersurface within
 %the product of spheres, and  $\mR / \R^+$ is the (skew) conormal bundle of this hypersurface. 
 \end{theorem}
 
 {\sc Proof.}   We first  check  that $\mR$ does not intersect the zero section.
 If $Q = 0$ then $0 =q_1 + t v_-$ which is only possible if either $q_1 = 0$ or
 $v_- \in L_1$ with $v_- = -q_1/t$.  The latter is impossible since this would imply that the whole incoming line
 $\ell_- \subset L_1$.  If the itinerary has length 2, the former is also possible, since
 if $q_1 = 0$, then $q_1 \in L_2$ as well, which is excluded by our definition of belonging to $\mB$.   Thus the quotient
 $\mR/ \R^+$ is a well-defined submanifold of $\bP T^* (S(E_v) \times (S(E_v) )$.
 
 Next we check the Legendrian condition and at the same time work out the contact form.
 Write $D = {\vec Q_-} \dd{}{\vec Q_-} + {\vec Q_-} \dd{}{\vec Q_-}$ for the Euler vector field, this being the vector field  whose flow is  dilation 
 (with $\lambda = e^t$ if $t$ is the flow parameter).  Let  $\Omega = \omega_- - \omega_+$
 be the symplectic form with respect to which $\mR$ is Lagrangian.   A standard construction  from contact and symplectic geometry
 suggests forming the one-form $\Theta = i_ D  \Omega$ which a direct  computation shows that 
 $\Theta$ is the one-form stated in the theorem.  Since $\mR$ is scale invariant, $D$ is tangent to $\mR$
 and consequently $\Theta (v) = \Omega(D, v) = 0$ for any other vector $v$ tangent to $\mR$.
 This proves that $\mR/\R^+$ is Legendrian relative to $\Theta$.  
 
 Finally, if $\mR/\R^+$ is transverse
 to the fibers of the fibration, then $\pi_- \times \pi_+$ maps it locally diffeomorphically onto a hypersurface.
 The projection and $\mR/\R^+$ are both  algebraic so if the Legendrian submanifold  is transverse at one point it is transverse
 at almost every point  and the image of each component is a singular hypersurface. 
$\Box$   
   %which are the subject of our main theorems  \ref{thm:main}, \ref{thm: main_unreduced} lie total collision ..
   
    \begin{remark}
    If $\mR/\R^+$ is nowhere transverse to the fiber then it is mapped to a subvariety of codimension greater than 1
    within the product of the spheres. 
    \end{remark}
    
 \begin{remark} We can summarize the discussion of this subsection as saying that  the map $q \mapsto (v_-, v_+)$
 which sends a billiard trajectory to its incoming and outgoing velocities is a ``Legendrian map''. 
 Arnol'd \cite[Appendix 16, p. 487]{Arnold} calls the restriction of a Legendrian fibration to a Legendrian submanifold
 a ``Legendrian map'' and its image a ``front'' as in ``wave-front.
 So the  ``scattering map''   $\pi_- \times \pi_+$  restricted to $\mR/ \R^+$ is
 a Legendrian map and its image, the `scattering front' will be an interesting singular hypersurface
 within the product of the incoming and outgoing velocity spheres.  See   subsection \ref{s:scatteringsurface}
below.   
  \end{remark}
  
%  \eject

\section{Examples}

\subsection{Origami unfoldings}
\label{s:examples}

\begin{figure}
\begin{center}
\vspace{-2mm}
\includegraphics[height=6cm]{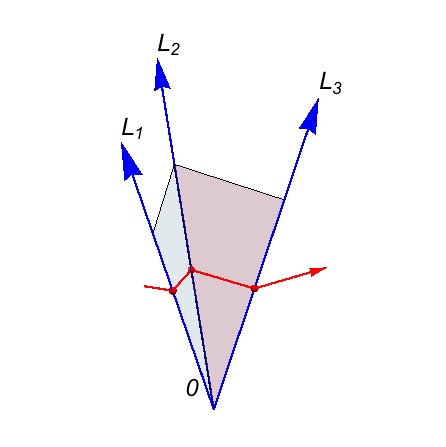}
\hfill
\includegraphics[height=6cm]{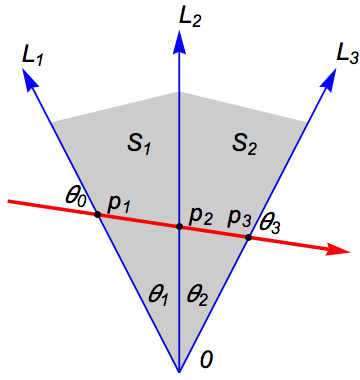}
\vspace{-5mm}
\end{center}
\caption{Left: Origami in $E$. Right: Origami unfolded}
\label{fig:origami}
\end{figure}

Suppose that   $\mL$
consists of   lines, so that $d = \dim(E) -1$.   Let $q \in \mB(L_1 L_2 \ldots L_k)$
be a  trajectory.  Then    each $q_i \in L_i$
is nonzero, for otherwise $q_i \in L_{i + 1}$ which is forbidden.    Assume $k > 1$. Edge   $q_i q_{i+1}$ of our $k+1$-gon $q$
joins the rays $\overrightarrow{0 q_i}$ and  $\overrightarrow{0 q_{i+1}}$ and hence lies in the plane $P_i$ spanned by $0, q_i, q_{i+1}$.
Within this plane the edge lies within the   sector
$S_i$ bounded by these two rays.    (By a ``sector'' we mean a planar convex region 
bounded by two rays.) Let  $\theta_i = {\rm angle}(q_i 0 q_{i+1})$ denote the opening angle of this
sector.  
Thus the interior part $q_1 q_2 ... q_k$ of our 
billiard trajectory lies on a polygonal cone within $E$ whose faces
are the sectors  $S_1 \ldots S_{k-1}$ glued together along the rays $0q_i \subset L_i$.  
We can ``unfold''  this cone onto a fixed plane, thus forming a big sector which is
made   of  congruent copies of our sectors $S_1, S_2, \ldots, S_{k-1}$  joined along their shared rays;
see figure \ref{fig:origami}. 

The opening angle of this big
developed sector is  
$$\beta = \theta_1 + \theta _2 + \ldots \theta_{k-1}.$$ 
Our billiard trajectory unfolds onto this developing plane as well.    
The billiard condition (1) is precisely that this unfolded trajectory is a straight line segment on this developing plane.
To reiterate,\\ 
{\it The billiard  segment $q_1 \ldots q_k$ becomes a   straight line segment  drawn 
on our big sector which is the  flattened polyhedral  cone! }

\begin{cor} \label{cor1}
 If $\beta \ge \pi$ the alleged billiard trajectory does not exist.
\end{cor}

\vskip .2 cm  
\noindent
{\sc Proof.} For if  a straight line  enters in to one part of a sector through one ray boundary and leaves through
the other ray boundary then the opening angle of the sector cannot be greater than $\pi$.
Said differently, the developed sector will not be convex if $\beta > \pi$. 
\hfill $\Box$\\[2mm]
Set  $\beta_i = {\rm angle}(L_i, L_{i+1})$ so that  $\beta_i$ is the minimum of $\theta_i$
and $\pi -\theta_i$.  Thus $\beta \ge \Sigma \beta_i$.

\begin{cor}  \label{cor2}
Set $\theta_{\min} = \min  ({\rm angle} (L, M))$, the minimum taken over all 
$L,  M \in \mL$, $L \ne M$.  There are no itineraries
of length greater than $1+ \lfloor \pi/\theta_{\min} \rfloor $.
\end{cor}

\vskip .2 cm  
\noindent
{\sc Proof.} Indeed since $\beta_i \ge \theta_{\min}$ we have that
$\beta \ge  (k-1) \theta_{\min}$ so that $\pi \ge (k-1) \theta_{\min}$ and
thus the number of intersection $k$ satisfies $(\pi/\theta_{\min})  \ge k -1$ . 
\hfill $\Box$\\[2mm]
Projecting the incoming and outgoing velocities onto our developing plane
we get information on their angles from the unfolded figure.

\begin{cor}  \label{cor3}
Consider the angle $\theta_0$ between the incoming ray (direction $v_A$) 
and the line $L_1$, the angle oriented so as to be the angle between 
$q_1$ and $-v_A$.  
Similarly consider the angle $\theta_k$ between
the final collision line $L_k$ and the outgoing ray (direction $v_B$),  
that angle oriented so as to be between
the vector $q_k$ and vector $v_B$.  
Then: 
$$\theta_{0} + \beta + \theta_k = \pi ; \qquad{ i.e.}\quad
  \Sigma_{i =0} ^k  \theta_i  = \pi . $$
\end{cor}

\vskip .2 cm  
\noindent
{\sc Proof.} 
Indeed, add on open planar sectors  with the plane $A L_1$ and $L_k B$ to the 
polygonal figure
described above, and flatten it.  Our billiard trajectory is a straight line on the resulting plane
and the angle sum, simply the opening angle of a line, is $\pi$.\hfill $\Box$\\[2mm]
We can now give a precise description of the Lagrangian relation 
$\tilde \mR (L_1 \ldots L_k)$.
For the relation to be nonempty we require $\Sigma \beta_i < \pi$.  
For each $i$ we consider two possible angles, $\beta_i$ and $\pi - \beta_i$.  
In all then, we have a collection of   $2^{k-1}$  angles $\theta_i$, each $\theta_i$ being either $\beta_i$
or $\pi - \beta_i$.  Among all these angle selections $\theta_1, \theta_2 , \ldots , \theta_{k-1}$ we  only consider  those
for which $\beta := \Sigma \theta_i < \pi$.   Now fix such a selection.  If our incoming line hits $L_1$ at an angle $\theta_0$, as defined in the corollary,
then our outgoing line must leave at angle $\theta_k = \pi - \beta - \theta_1$. 
Let $q_1, q_k$ be the points where the incoming line hits $L_1$ and where the outgoing line
leaves $L_k$ and let $v_1, v_k$ be the corresponding velocities.
We  identify the lines $(\ell_-, \ell_+) \in \mR$  according to their
boundary conditions $(q_1, v_1)$, $(q_k,  v_k)$. 
We have shown that
$$\theta_0 + \beta + \theta_k = \pi$$
Referring to figure \ref{fig:origami}
we have, by the law of sines:
$$|q_1|/\sin(\theta_k) = |q_k|/\sin(\theta_0).$$
These two relations, together with the specification of 
$\beta$, then determine our Lagrangian relation.

\vskip .2 cm  
\noindent
{\sc Generalizations.}

Take  now an arbitrary collection $\mL$ of subspaces of the same dimension. 
 
Let $A q_1 q_2 \ldots q_k B$ be any billiard trajectory.
Then   Cor.\ \ref{cor1} and Cor.\  \ref{cor3} hold, with the angles $\theta_i$ now being ${\rm angle}(q_i, q_{i+1})$.  Indeed, just take the $L_i$ to be the rays $\overrightarrow{0q_i}$ and proceed as before! 

Cor.\  \ref{cor2} generalizes to the case in which $\mL$ is comprised of higher-dimensional subspaces, instead of lines, {\it provided
these subspaces enjoy the property that the  intersection of any two of them    is zero}.
For two subspaces $L, M$  whose intersection is zero,  
the notion of the minimal angle between them makes sense:
$\theta (L, M) = \min_{v,w}{\rm angle} (v, w)$ with the minimum taken over all nonzero pairs
$v \in L,  w \in M$.    Now, in this setting we define 
$\theta_{\min}$ as above.  Corollary  \ref{cor2} holds exactly as stated. 

%\eject

\subsection{A Scattering Surface}
\label{s:scatteringsurface}

This  example  illustrates the use of symmetry and scaling (section \ref{s:symmetries}) and  the complexity of the scattering relation  
even for apparently simple itineraries.  Working out this example inspired  the discovery of theorem \ref{Legendrian}. 

Consider 3-body billiards for  three equal masses moving in the plane
$\R^2 = \C$.   Write $q = (q_1, q_2, q_3) \in \C^3$ for the positions of the masses
and $v = (v_1, v_2, v_3) \in \C^3$ for velocities.
We  set the linear mometum equal to zero:
$v_1 + v_2  + v_3 = 0$
and assume that the center of mass is also zero.
In this way, the underlying Euclidean space $E$
 becomes   real 4 dimensional (or complex 2-dimensional) linear subspace of   $\C \times \C \times \C$.
 Our collision subspaces are the $\Delta_{ij}$ intersected with this $E$. 
  \begin{figure}
\begin{center}
\includegraphics[height=8cm]{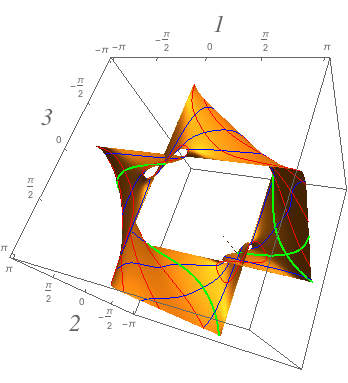}
\vspace{-5mm}
\end{center}
\caption{A slice $\pi_+ (\pi_- ^{-1}(v))$ of the scattering surface.}
\label{fig:scatteringDirns}
\end{figure}

Fix the itinerary to be $\Delta_{12} \Delta_{13} $:   first 1 and 2 collide, then 1 and 3.
We will explore a small  part of the corresponding translation-reduced scattering relation 
$\mR =\mR(\Delta_{12} \Delta_{23})$.  The space of lines in $\R^4$ being 6 dimensional,  $\mR$ forms a 6-dimensional Lagrangian relation, 
so a 6-dimensional submanifold of   $T^* S^3 \times T^* S^3$.We will fix the incoming direction $v_-$  of the incoming ray $\ell_-$ but we leave the ``impact parameter'
$Q_-$ free. For specificity, let us  fix the incoming direction $v_- \in E_v$
by supposing that 
the  three masses come in from infinity  with their directions equally spaced to
form the   vertices of an equilateral triangle:
$(v_{1-}, v_{2-}, v_{3-}) = (1 ,  exp (2 \pi i/3), exp(4 \pi i / 3)$.  
(Take the equal masses to be $m_1 = m_2 = m_3 = 1/3$ so that  $v_{-}$ is unit length.)
We computed all possible outgoing velocities $v_+ = (v_1 ^+, v_2 ^+, v_3 ^+)$.  The results are depicted in \ref{fig:scatteringDirns}.
 These outgoing velocities form a 2-dimensional surface within
the 3-sphere $S(E_v)$ of all possible unit length velocities in $E$. We have 
coordinatized this surface by projecting $v_+$ to its three component vectors $v_i ^+ \in \C$ and  then plotting the argument of that complex number.  

As described in theorem \ref{Legendrian}, the quotient of $\mR$ by scaling forms a 5-dimensional   Legendrian submanifold $\mR/ \R^+$ inside the   projectivized cotangent bundle
of the product of our incoming and outgoing velocity spheres.  
The Legendrian map $\pi_- \times \pi_+$ of theorem \ref{Legendrian} takes the relation
onto a 5 dimensional   hypersurface (probably with singularities)  within the product of the incoming and outgoing velocity 3-spheres
By freezing the value of $v_-$ we have depicted  in figure \ref{fig:scatteringDirns} a single two-dimensional  `slice'  of this hypersurface,
namely the surface  $\pi_+ (\pi_- ^{-1} (v_-))$. 

Here are  details of the computation leading to the figure. 
After the 1st collision of 1 with 2, the 3rd particle's velocity is unchanged.
Write $v^m = (v_1 ^m , v_2 ^m, v_3 ^m)$  for this intermediary velocity, between $\Delta_{12}$ and $\Delta_{23}$.
Then  $v_3 ^m = v_3 ^-$ and    there
is a vector $w \in \R^2$ such that 
$v_1^m = \eh(v_1^-+v_2^-)+ w$, 
$v_2^m = \eh(v_1^-+v_2^-) - w$ with $|w| =\eh |v_1^- - v_2^-|
= \frac{3}{2}$ and the direction of $w$ arbitrary. 

After collision of particle number 1 and 3,
we get our final velocity $v^+ = (v_1 ^+, v_2 ^+, v_3 ^+)$
with   $v_2^+=v_2^m$,
$v_1^+=\eh(v_1^m+v_3^m)+u$, $v_3^+= \eh(v_1^m+v_3^m)-u$  
with $|u| =\eh |v_3^m - v_1^m|$.  Use $v_1 + v_2 + v_3 = 0$
so that $v_1 ^m = -\eh v_3^- + w$ to rewrite
$v_3^m - v_1^m = \frac{3}{2} v_3 ^- -w$ so that $|u| = | \frac{3}{2} v_3 ^- -w|$.

\section{Open Problems}

\subsection{ On the closure of the Lagrangian relations.}  
$ $

\vskip .4 cm  

\noindent
{\sc Question 1.}   What are the closures of the Lagrangian relations of theorem \ref{main}? 

Recall that these relations are denoted $\mB(L_1 \ldots L_k)$ where $L_1 L_2 \ldots L_k$ is the itinerary.

\begin{example} Let us suppose the codimension $d > 1$ and 
that  $\mB(L_1 L_2) \ne \emptyset$.  Then it must be that $\mB(L_1) \ne \emptyset$ and moreover
${\rm cl}(\mB(L_1 )) \cap \mB(L_1 L_2) \ne \emptyset$.  
For suppose that $q \in \mB (L_1 L_2)$. Then $q$
has an edge $q_1 q_2$ joining $L_1$ to $L_2$.  We can perturb the endpoint $q_2$
slightly, off into $E$, and insure that the resulting ray $\overrightarrow{q_1 q_2}$ never intersect $C$ again.
\end{example}

\vskip .2 cm  
\noindent

{\sc Question 2.} What algebraic or combinatorial relationships hold between our  Lagrangian relations?   

Lagrangian relations are built to be composed.   (See for example \cite{GS} on composing linear Lagrangian relations.) How and when can we compose our Lagrangian relations?
Concatenation of polygonal paths
{\it suggests} that  their should be some type of composition
law 
$$\mB(L_1 L_2 \ldots L_k)  \times  \mB (L_{k_1} L_{k+2} L_{k+s}) \to ^{?} \mB(L_1 L_2 \ldots L_k L_{k+1} \ldots L_{k+s}) $$
This ``law''  is  nonsense if taken literally.  Indeed it is doomed to failure by   the background theorem,
theorem \ref{sss:FFT} which implies that   concatenations between relations
cannot be arbitrarily long for their target is then empty.        

It seems  there does exist, however, some kind of ``decomposition''.  
  Write $I =  L_1 L_2 \ldots L_1 L_2 \ldots L_k,  J = L_{k+1} \ldots L_{k+s}$
for two itineraries.  Suppose
we have a path  $q \in \mB (I J )$.
Moving forward in time along $q$, at each collision point $q_i \in L_i$ we have a continuous choice of new outgoing directions.  In particular,
at the $k$th step we could make this choice so that the new outgoing ray never intersects $C$
again.  In this way we would achieve, by perturbing $q$ at the $k$th collision,  a $\tilde q \in \mB (L_1 \ldots L_k)$.  
Thus there appears to be a well-defined map: 
$\pi:  \mB (IJ )  \times U \to \mB (I)$
where $U$ are   `perturbation parameters'' describing how we perturbed the final outgoing ray from $L_k$ so as to sail off to infinity. 
Presumably $U \subset S^{d-1}$.  Viewing this same  perturbation `process'' backwards in time, we could vary the incoming direction to $L_k$
at $q_k$ to arrive at a $\tilde q_-  \in \mB (J)$, and so arrive at a `decomposition map'' $\mB (I J ) \times \tilde U \to \mB(I) \times \mB (J)$. 

SUBQUESTION:  Is there a well-defined decomposition ``morphism'' $\mB (IJ )   \to \mB (I) \times \mB(J)$ ?
%for appropriate pairs $I, J$, maybe with additional parameters required to have a good map?

\vskip .2 cm  
\noindent

{\sc QUESTION 3.}  
List all the possible itineraries $I$  with nonempty realizations
$\mB(I)$?

We know by the background  theorem \ref{sss:FFT}  that this list is finite.
%Our interest is in the $N$-body case.  
This last question seems to be the   simplest, hardest question we have asked so far. 
  
%\section{acknowledgement}

\begin{acknowledgement}RM would like to thank Alan Weinstein for early interest in this work, Mikhail Kapovich
for discussions of CAT(0), 
and Bruno Nachtergaele for a crucial question leading to the section on conservation laws.
RM thankfully acknowledges support by NSF grant DMS-20030177.
\end{acknowledgement}

\end{document}